\documentclass[11pt,british]{amsart}

\newif\ifpdf
\ifx\pdfoutput\undefined\pdffalse\else\pdfoutput=1\pdftrue\fi\newif\pdf
\ifpdf\relax\else
\usepackage[T1]{fontenc}
\fi

\usepackage{babel,eucal,url,amssymb,stmaryrd,booktabs,%
enumerate,amscd,paralist}

\theoremstyle{plain}
\newtheorem{thm}{Theorem}[section]

\newtheorem{pro}[thm]{Proposition}
\newtheorem{co}[thm]{Corollary}

\theoremstyle{definition}

\theoremstyle{remark}
\newtheorem{rem}[thm]{Remark}
\newtheorem{example}[thm]{Example}
\newtheorem*{ack}{Acknowledgements}

\newcommand{\Lie}[1]{\operatorname{\textsl{#1}}}
\newcommand{\lie}[1]{\operatorname{\mathfrak{#1}}}
\newcommand{\GL}{\Lie{GL}}

\newcommand{\SO}{\Lie{SO}}
\newcommand{\so}{\lie{so}}

\newcommand{\Spin}{\Lie{Spin}}
\newcommand{\SU}{\Lie{SU}}

\newcommand{\Un}{\Lie{U}}
\newcommand{\un}{\lie{u}}
\newcommand{\Gtwo}{\ifmmode{{\rm G}_2}\else{${\rm G}_2$}\fi}

\newcommand{\sumcic}{\mathop{\mbox{\large {$\mathfrak S$}\vrule width 0pt
depth 2pt}}}

\def\sideremark#1{\ifvmode\leavevmode\fi\vadjust{\vbox to0pt{\vss
 \hbox to 0pt{\hskip\hsize\hskip1em
 \vbox{\hsize2.5cm\tiny\raggedright\pretolerance10000
 \noindent #1\hfill}\hss}\vbox to8pt{\vfil}\vss}}}%

\date{\today}
\begin{document}

\title[$\SU(3)$-structures on submanifolds of a SPIN(7)-manifold]%
{$\SU(3)$-structures on submanifolds of a SPIN(7)-manifold}

\author{Stefan Ivanov}
\address[S.Ivanov]{University of Sofia "St. Kl. Ohridski"\\
Faculty of Mathematics and Informatics\\
Blvd. James Bourchier 5\\
1164 Sofia, Bulgaria} \email{ivanovsp@fmi.uni-sofia.bg}

\author{Francisco Mart\'\i n Cabrera}
\address[F.Mart\'\i n Cabrera]{Department of Fundamental Mathematics\\
  University of La Laguna\\ 38200 La Laguna, Tenerife, Spain}
\email{fmartin@ull.es}

\begin{abstract}
Local $\Lie{SU}(3)$-structures on an oriented submanifold of
$\Lie{Spin}(7)$-manifold are determined and their types are
characterized in terms of the shape operator and the type of
$\Spin(7)$-structure.  An application to Bryant \cite{Br3} and
Calabi \cite{Cal} examples is given. It is shown that the product
of a Cayley plane and a minimal surface lying in a
four-dimensional orthogonal Cayley plane with the induced complex
structure from the octonions
 described by Bryant in \cite{Br3} admits a
holomorphic local  complex volume form exactly when it lies in a
three-plane, i.e. it coincides with the example constructed by
Calabi in \cite{Cal}. In this case the holomorphic $(3,0)$-form is
parallel with respect to the unique Hermitian connection with
totally skew-symmetric torsion.



\end{abstract}

\maketitle \setcounter{tocdepth}{2} \tableofcontents

\section{Introduction}\label{introduction}
A $\Lie{Spin}(7)$-structure on an eight-dimensional manifold is by
definition a reduction of the structure group of the tangent
bundle to $\Lie{Spin}(7)$. An eight-dimensional manifold equipped
with a $\Lie{Spin}(7)$-structure is called
$\Lie{Spin}(7)$-manifold. Moreover, associated with a
$\Spin(7)$-structure, there exists a nowhere vanishing four-form
$\Phi$, called \emph{the fundamental form}, which determines  a
Riemannian metric $\langle \cdot  , \cdot \rangle$ and a volume
form due to the fact that $\Lie{Spin}(7)$ is the maximal compact
subgroup of $\Lie{SO}(8)$. Likewise, choosing a vector of unit
length   as unity, the tangent vector space on each point of a
$\Lie{Spin}(7)$-manifold can be identified with the octonian
algebra ${\mathbb O}$.

Decomposing the space $\{\nabla\Phi\}$  of  covariant derivatives
of $\Phi$ with respect to the Levi-Civita connection $\nabla$ into
$\Lie{Spin}(7)$-irreducible components, Fern{\'a}ndez \cite{F}
classified $\Lie{Spin}(7)$-manifolds  and obtained four classes,
namely, $\overline{W}_0$ (parallel), $\overline{W}_1$ (balanced),
$\overline{W}_2$ (locally conformal parallel) and the whole class
$\overline{W}$. Studying the obstruction of
 a $\Lie{Spin}(7)$-structure to be parallel, we  find (Theorem~\ref{spin7t})
an expression for the intrinsic torsion of a
 $\Lie{Spin}(7)$-structure in terms of
the exterior derivative $d \Phi$ which explicitly expresses
 $\nabla\Phi$ in terms of $d\Phi$.
Note, that a formula
 of $\nabla\Phi$ in terms of $d \Phi$ was given in \cite{I1}.
 The existence of such an explicit formula is an implicit
 consequence
of the fact, noted by Bryant \cite{Br} (see \cite{F,Sal}), that
the Riemannian holonomy group of a $\Lie{Spin}(7)$-manifold is
contained in $\Lie{Spin}(7)$ iff the form $\Phi$ is closed.

If $M^6$ is an orientable six-dimensional submanifold of a
$\Lie{Spin}(7)$-manifold $(M^8,\Phi,\langle \cdot , \cdot \rangle
)$, Gray \cite{Gr} showed that there is on $M^6$ an almost
Hermitian structure ($\Lie{U}(3)$-structure)  naturally induced
from the $\Lie{Spin}(7)$-structure on $M^8$. When $M^8$ is a
parallel $\Lie{Spin}(7)$-manifold, Gray  derived conditions in
terms
 of the shape operator of $M^6$ characterizing types of almost Hermitian structure on
 $M^6$.

In the present paper, we define local $\Lie{SU}(3)$-structures on
$M^6$ inherited from the $\Lie{Spin}(7)$- structure on $M^8$. Note
that in general there is not a  global $\SU(3)$-structure on $M^6$
induced from the $\Lie{Spin}(7)$ structure on $M^8$,  since the
stabilizer of an oriented two-plane in $\Lie{Spin}(7)$ is the
group $\Lie{U}(3)$ \cite{Br3}. We show the existence of local
complex volume forms naturally induced from the fundamental
four-form $\Phi$ and the choice of a local oriented orthonormal
frame $N_1$, $N_2$ of the normal bundle of $M^6$. We present
relations between the $\Lie{Spin}(7)$-structure on the ambient
manifold $M^8$, the induced local $\Lie{SU}(3)$-structure and the
shape operator on $M^6$ (Proposition \ref{rraa}). Consequently, we
characterize the types of the local $\Lie{SU}(3)$-structures on
$M^6$ in terms of the fundamental four-form $\Phi$ and the shape
operator  (Theorem \ref{wzerogeneral}, Theorem \ref{wonegeneral},
Theorem \ref{wtwogeneral}).  In particular, we recover  Gray's
results in \cite{Gr} in an alternative way.

In  Section~\ref{I1} we study the problem  when there exists a
closed local  $\Lie{SU}(3)$-structure on $M^6\subset M^8$, which
in particular, implies that the almost complex structure is
integrable due to the considerations in \cite{hitc}. We focus our
attention to the case $M^8={\mathbb O}$ studied in detail by
Bryant in \cite{Br3}. In this case (even more general, when the
$\Lie{Spin}(7)$-structure of the ambient manifold is parallel),
Gray \cite{Gr} showed that the Lee form of the submanifold  is
always zero. When the almost complex structure is integrable, then
it is balanced (type $\mathcal{W}_3$) and  the submanifold is
necessarily minimal. The properties of submanifolds with balanced
Hermitian structure are investigated by Bryant in \cite{Br3}. He
shows that if $M^6\subset{\mathbb O}$ inherits complex and
non-K{\"a}hler structure, then  $M^6$ is foliated by four-planes
in ${\mathbb O}$ in a unique way, he calls this foliation
asymptotic ruling. He  obtains that if the asymptotic ruling is
parallel,  then $M^6$ is a product of a fixed associative
four-plane $Q^4$ in ${\mathbb O}$ with a minimal surface in the
orthogonal four-plane. He shows that the Calabi examples,
described in \cite{Cal}, are exactly those complex $M^6$ with
parallel asymptotic ruling contained in $Im{\mathbb O}\subset
{\mathbb O}$, i.e. the minimal surface lies in an associative
three-plane in $Im {\mathbb O}$.

We investigate  when there exists a local holomorphic $\Lie{SU}(3)$-structures
 in the case of parallel asymptotic ruling. We show
that there exists a holomorphic local $\Lie{SU}(3)$-structure on
$M^6$ exactly when the minimal surface lies in a three-plane
(Theorem~\ref{cal}). We also prove that the corresponding Bismut
connection (the unique Hermitian connection with totally
skew-symmetric torsion) preserves the holomorphic volume form
having holonomy contained in $\Lie{SU}(3)$. Therefore, the
structure is  Calabi-Yau with torsion  (CYT). CYT structures are
attractive in heterotic string theory as a possible solution to
the heterotic string model proposed by Str\"ominger \cite{stro}.
Consequently, we derive that the compact complex non-K{\"a}hler
six-manifold with vanishing first Chern class constructed by
Calabi in \cite[Theorem 7]{Cal} has holomorphically trivial
canonical bundle and the $\Lie{SU}(3)$-structure constructed by
Calabi is a CYT-structure (Theorem~\ref{cal1}).

Recently,  Bryant discussed in \cite{Br1} a generalization of the
notion of holomorphic vector bundles on complex manifold to the
almost complex case and,  consequently, a generalization of the
notion of Hermitian-Yang-Mills connection. He referred the class
of almost complex six-manifolds admitting such non-trivial bundles
as quasi-integrable. An important subclass is the class of strict
quasi-integrable structures which is defined as quasi-integrable
structures with nowhere vanishing Nijenhuis tensor. He introduced
the notion of quasi-integrable $\Lie{U}(3)$-structure, pointing
out that this class of almost Hermitian six-manifold coincides
with the class $\mathcal W_1\oplus \mathcal W_3\oplus \mathcal
W_4$ according to Gray-Hervella classification \cite{GrH}, i.e.
the class where the Nijenhuis tensor is totally skew-symmetric.
The case of nearly K\"ahler structures is also investigated in
details in \cite{Br1}. Following our approach, in
Example~\ref{exam1}, we describe a strict quasi-integrable
non-nearly K\"ahler $\Lie{SU}(3)$-structures on $S^3\times S^3$
compatible with the standard product metric on $S^3\times S^3$.
Four of these structures are half-flat in the class $\mathcal
W_1\oplus \mathcal W_3$. These four structures are also
left-invariant on the group $SU(2)\times SU(2)\cong S^3\times
S^3$.
\begin{rem}
We note that another compact example of strict quasi-integrable
non-nearly K\"ahler half-flat $\Lie{SU}(3)$-structure of type
$\mathcal W_1\oplus \mathcal W_3$ tensor on nil-manifold has been
constructed in \cite{II}, Section~6.2.
\end{rem}

\begin{ack}  The second author is supported by a grant from MEC (Spain), project
  MTM2004-2644. Stefan Ivanov  thanks University of
  La Laguna for the very  kind hospitality during the initial stages of this
  work. Likewise, he also thanks Thomas Friedrich, Ilka Agricola and the
junior research group (VolkswagenStiftung) around I. Agricola at
Humboldt University in Berlin where a part of this work was
  done.
\end{ack}

\section{General properties of $\SU(3)$ and $\Spin(7)$-structures}\label{sec:preliminaries}
 In this section we recall necessary properties of $\Lie{SU}(3)$
and $\Lie{Spin}(7)$-structures.

First we recall some notions relative to $G$-structures, where
$\Lie{G}$ is a subgroup of the linear group $\Lie{GL}(m , \mathbb
R)$. If $M$ possesses a $G$-structure, then there always exists a
$\Lie{G}$-connection defined on $M$. Moreover, if $(M^m ,\langle
\cdot , \cdot \rangle)$ is an orientable $m$-dimensional
Riemannian manifold with associated Levi-Civita connection
$\nabla$ and $G$ is  a closed and connected subgroup of $\SO(m)$,
then there exists a unique metric $G$-connection $\nabla^G$ such
that $\xi^G_x = \nabla^G_x - \nabla_x$ takes its values in
$\lie{g}^{\perp}$, where $\lie{g}^{\perp}$ denotes the orthogonal
complement in $\so(m)$ of the Lie algebra $\lie{g}$ of $G$
\cite{Sal,CleytonSwann:torsion}. The tensor $\xi^G$ is called the
{\it intrinsic torsion}  of the $G$-structure and $\nabla^G$ is
referred as  the {\it minimal $G$-connection}.

\subsection{$\SU(3)$-structures}  \label{su3structure}
Here we give a brief summary of the properties of
$\Lie{SU}(3)$-structures on six-dimensional manifolds which are
also called \emph{special almost Hermitian six-manifolds}.
 For more detailed and exhaustive information see \cite{CS,Cabrera:special}.

An {\it almost Hermitian} manifold is a $2n$-dimensional manifold
$M$ with a $\Un(n)$-structure. This means that $M$ is equipped
with a Riemannian metric $\langle \cdot , \cdot \rangle$
 and an
orthogonal almost complex structure $J$. Each fibre $T_m M$ of the
tangent bundle can be consider as complex vector space letting $i
x = Jx$. The K{\"a}hler form $\omega$ is defined by $\omega (x,
y)= \langle x, Jy \rangle$.

{\bf Convention.} For a $(0,s)$-tensor $B$, we write
\begin{equation*} \label{ecuacionesb}
  \begin{array}{l}
    J_{(j)}B(X_1, \dots, X_j, \dots , X_s) = - B(X_1, \dots , JX_j, \dots ,
    X_s),\\[2mm]
    J B(X_1,\dots,X_s) = (-1)^sB(JX_1,\dots,JX_s).
  \end{array}
\end{equation*}
The Lee form $\theta$ of an almost Hermitian structure is defined
by $\theta=Jd^*\omega$, where $d^*$ denotes the codifferential.
Also we will consider the natural extension of the metric $\langle
\cdot , \cdot \rangle$ to $\Lambda^p T^* M$ given by
$$
\langle \alpha , \beta \rangle = \tfrac{1}{p!} \sum_{i_1,\dots
,i_p=1}^{2n} \alpha ( e_{i_1} , \dots , e_{i_p} ) \beta ( e_{i_1}
, \dots , e_{i_p} ),
$$
where $\{e_1, \dots, e_{2n}\}$ is an orthonormal basis for
vectors.

A {\it special almost Hermitian} manifold is a $2n$-dimensional
manifold $M$ with an $\SU(n)$-structure. This means that $(M,
\langle \cdot, \cdot \rangle, J)$ is an almost Hermitian manifold
equipped with a complex volume form $\Psi = \Psi_+ + i \Psi_-$,
i.e. $\Psi$ is an $(n,0)$-form
 such
that $\langle  \Psi , \overline{\Psi} \rangle = 1$, where $\langle
\cdot , \cdot \rangle$ denotes the natural extension of the metric
on (complex) forms and $\overline{\Psi}$ is the conjugated
$(0,n)$-form.  Note that $J_{(j)} \Psi_+ = \Psi_-$. \vspace{2mm}

In general, an almost Hermitian manifold admits a linear
connection preserving the almost Hermitian structure and having
totally skew-symmetric torsion exactly when the Nijenhuis tensor
is totally skew-symmetric (the class $\mathcal W_1\oplus \mathcal
W_3\oplus \mathcal W_4$  in the Gray-Hervella classification
\cite{GrH}). Moreover, such a connection is unique
 \cite{FI,Gau}. If  the almost complex structure is integrable,  then
 this connection is referred as the Bismut connection.
It was  used by Bismut \cite{bis} to derive a local index formula
for Hermitian non-K\"ahler manifolds. When the Bismut connection
preserves a given  $\Lie{SU}(n)$-structure, i.e. it has holonomy
contained in $\Lie{SU}(n)$, then the manifold is called sometimes
Calabi-Yau manifold with torsion (CYT) and appears as a possible
geometry in heterotic string model due to the work of Str\"ominger
\cite{stro} (see e.g. \cite{BBE,Car,GMPW,GMW,GPap,GP,GLMW,GM,II}
and references therein).

In the following, we  consider special almost Hermitian
six-manifold, i.e. a six-dimensional smooth manifold endowed with
an $\Lie{SU}(3)$-structure. We denote the corresponding Lee form
by $\theta^6$. Let
$$e_{1\mathbb C}=e_1+iJe_1,\quad e_{2\mathbb C}=e_2+iJe_2,\quad e_{3\mathbb C}=e_3+iJe_3$$
be a unitary basis
  such that $  \Psi (e_{1\mathbb C}, e_{2\mathbb
C}, e_{3\mathbb C})=1$,  i.e. $ \Psi_+ (e_1, e_2 , e_3) = 1$,
 $\Psi_- (e_1 , e_2 , e_3)=0$. The real orthonormal basis for
vectors $e_1$,$e_2$,$e_3$,$Je_1$,$Je_2$,$Je_3$ is said to be {\it
adapted to the $\SU(3)$-structure}. By means of such an adapted
basis, the K\"ahler form $\omega$ and the three-forms
$\Psi,\Psi_+$ and $\Psi_-$ are given by
\begin{eqnarray*}
\omega&=&-e_1\wedge Je_1-e_2\wedge Je_2-e_3\wedge Je_3,\\
\Psi & = & e_{1\mathbb C} \wedge  e_{2\mathbb C} \wedge  e_{3\mathbb C},\\
\Psi_+ & = & e_1 \wedge e_2 \wedge e_3 - Je_1 \wedge Je_2 \wedge
e_3
-Je_1 \wedge e_2 \wedge Je_3 -e_1 \wedge Je_2 \wedge Je_3, \\
\Psi_- & = & - Je_1 \wedge Je_2 \wedge Je_3 + Je_1 \wedge e_2
\wedge
e_3 + e_1 \wedge Je_2 \wedge e_3 + e_1 \wedge e_2 \wedge Je_3.
\end{eqnarray*}
Here and further we freely identify vector field with the dual
one-form via the metric.

It is straightforward to check $
 \omega^3: = \omega \wedge \omega \wedge  \omega =
 6 \, e_1 \wedge e_2 \wedge e_3 \wedge Je_1 \wedge Je_2 \wedge Je_3.
$ If we fix the real volume form $Vol$
 such that $6 Vol = \omega^3$, we have the relations \cite{CS,Cabrera:special}
\begin{eqnarray} \label{volumenes}
\Psi_+ \wedge \omega & =& \Psi_- \wedge \omega = 0;\\\nonumber
 \Psi_+ \wedge \Psi_-& =& -   4 \, Vol, \qquad
\Psi_+ \wedge \Psi_+ = \Psi_- \wedge \Psi_- = 0;\\\nonumber
x \wedge \Psi_+ =  Jx \wedge \Psi_- &=&  - (Jx \lrcorner \Psi_+) \wedge
\omega, \qquad  x \lrcorner \Psi_+ = Jx \lrcorner \Psi_-, \quad x\in T_mM,\nonumber
\end{eqnarray}
where $\lrcorner$ denotes the interior product of vectors and
forms.

Note that, defined on $M$, there are two Hodge star operators
associated with the  volume forms $Vol$ and  $\Psi$. Relative to
the real Hodge star operator $\ast$,  for any one-form
$\mu\in\Lambda^1M$, we have  the relations
\begin{gather} \label{estrella1}
\ast \left( \ast ( \mu \wedge \Psi_+ ) \wedge \Psi_+ \right) =
\ast \left( \ast ( \mu \wedge \Psi_-) \wedge \Psi_- \right)  = -2
\mu, \\
 \ast \left( \ast ( \mu \wedge \Psi_- ) \wedge
\Psi_+ \right) = - \ast \left( \ast ( \mu \wedge \Psi_+ ) \wedge
\Psi_- \right) =  2
J\mu. \label{estrella2}
\end{gather}

For $\Un(3)$-structures, the minimal $\Un(3)$-connection is given
by $ \nabla^{\Un(3)} = \nabla + \xi^{\Un(3)}$, where
\begin{equation} \label{torsion:xi}
\xi^{\Un(3)}_X Y = - \frac{1}{2} J \left( \nabla_X J \right) Y
\end{equation}
(see \cite{Falcitelli-FS:aH}). Since $\Un(3)$ stabilizes the
K{\"a}hler form $\omega$, it follows that $\nabla^{\Un(3)} \omega
= 0$. Then $\nabla \omega = - \xi^{\Un(3)} \omega \in T^* M
\otimes \un(3)^{\perp}$. Thus, one can identify the
$\Un(3)$-components of $\xi^{\Un(3)}$ with the $\Un(3)$-components
of $\nabla \omega$.

For $\SU(3)$-structures, we have the decomposition $\lie{so}(6)=
\lie{su}(3) \oplus \mathbb R \omega \oplus \un(3)^{\perp}$,
 i.e.
$\lie{su}(3)^{\perp} = \mathbb R \omega \oplus \un(3)^{\perp}$.
Therefore, the intrinsic $\SU(3)$-torsion $\eta + \xi^{\Un(3)}$ is
such that $\eta \in T^* M \otimes \mathbb R J \cong T^* M$
 and $\xi^{\Un(3)}$ is still determined by
Equation \eqref{torsion:xi}. The tensors $\omega$, $\Psi_+$ and
$\Psi_-$ are stabilized by the $\SU(3)$-action and therefore
$\nabla^{\SU(3)} \omega = 0$, $\nabla^{\SU(3)} \Psi_+ = 0$, $
\nabla^{\SU(3)} \Psi_- = 0$, where
$$
\nabla^{\SU(3)} = \nabla + \eta + \xi^{\Un(3)}
$$
is the minimal $\SU(3)$-connection. Since $\nabla^{\SU(3)}$ is
metric and $\eta \in T^* M \otimes \mathbb R J$, we have $ \langle
Y , \eta_X Z \rangle = (J \eta)(X) \omega(Y,Z)$, where $\eta$ on
the right hand side is considered to be  a one-form. Hence
 \begin{equation} \label{torsion:eta}
 \eta_X Y = J \eta (X) JY.
\end{equation}

One can check  $\eta \omega = 0$, then  from $\nabla^{\SU(3)}
\omega = 0$ one gets
$$
 \nabla \omega = - \xi^{\Un(3)} \omega \in T^* M
\otimes \lie{u}(3)^{\perp} = {\mathcal W}_1 \oplus {\mathcal
W}_2\oplus {\mathcal W}_3 \oplus {\mathcal W}_4,
$$
where the summands ${\mathcal W}_i$ are the  Gray-Hervella
$\Un(3)$-modules.
 There
is a further splitting of $T^* M \otimes \lie{u}(3)^{\perp}$ into
six  $\Lie{SU}(3)$-modules  discovered and first described by
Chiossi and Salamon in \cite{CS} (see also
\cite{Cabrera:special,Cabrera:horient}, for interpretation in
physics see \cite{Car,GLMW,GM}). We present below the necessary
for our considerations part of the description of the
$\Lie{SU}(3)$-modules following \cite{Cabrera:special}.

The spaces ${\mathcal W}_3$ and ${\mathcal W}_4$ are irreducible
also as $\SU(3)$-modules. However, ${\mathcal W}_1$ and ${\mathcal
W}_2$ admit the  decompositions $ {\mathcal W}_j = {\mathcal
W}_j^+ \oplus {\mathcal W}_j^-$, $j=1,2$, into irreducible
$\SU(3)$-components, where ${\mathcal W}_j^+$ (resp. ${\mathcal
W}_j^-$) includes those  elements $\beta \in {\mathcal W}_j
\subseteq T^*M \otimes \Lambda^2 T^* M$ such that the bilinear
form $r(\beta)$, defined by $2 r(\beta)(x,y) = \langle  x
\lrcorner \beta , y \lrcorner \Psi_+ \rangle$, is symmetric
 (resp. skew-symmetric).

On the other hand, we have  $$ \nabla \Psi_+ = - \eta \Psi_+ -
\xi^{\Un(3)} \Psi_+, \quad \nabla \Psi_- = - \eta \Psi_- -
\xi^{\Un(3)} \Psi_-,$$
 since $\Psi_{\pm}$ are $\nabla^{\SU(3)}$-parallel.

 Therefore, using   \eqref{torsion:xi}
and \eqref{torsion:eta}, we obtain the following expressions
\begin{equation*} \label{torsiones}
\begin{array}{lll}
\quad - \eta_X \Psi_+   =  - 3 J\eta(X) \Psi_-, & \quad &
-\xi^{\Un(3)}_X \Psi_+  = \displaystyle
\frac{1}{2}\sum_{j=1}^6((e_j \lrcorner \nabla_X \omega) \wedge
(e_j \lrcorner \Psi_-)), \\[2mm]
\quad - \eta_X \Psi_-   =    3 J\eta(X) \Psi_+,  & &
 -\xi^{\Un(3)}_X \Psi_-  = - \displaystyle \frac{1}{2}\sum_{j=1}^6(  (e_j \lrcorner \nabla_X \omega)
\wedge (e_j \lrcorner \Psi_+)),
\end{array}
\end{equation*}
where $\{ e_1 , \dots , e_6 \}$ is an orthonormal basis for
vectors.

 Denote $ {\mathcal W}^{\pm}_{5}
 = T^* M \otimes \Psi_{\pm}$.
 It is clear that $- \eta \Psi_+ \in {\mathcal W}^-_{5}$, $- \eta \Psi_- \in {\mathcal W}^+_{5}$.

Consider the two
$\SU(3)$-maps $
 \Xi_+ , \Xi_- \, : \, T^* M \otimes \lie{u}(3)^{\perp} \to T^* M \otimes \Lambda^3 T^* M $ defined by
$$\nabla_{\cdot} \omega \to \frac12 \,\sum_{j=1}^6 (e_j \lrcorner
\nabla_{\cdot} \omega) \wedge (e_j \lrcorner \Psi_-)), \quad \nabla_{\cdot}
\omega \to  - \frac12 \,\sum_{j=1}^6(
 (e_j \lrcorner \nabla_{\cdot} \omega) \wedge (e_j \lrcorner \Psi_+)),$$ respectively. It turns out that
the $SU(3)$-maps $\Xi_+$ and $\Xi_-$ are injective and
\begin{equation*}
\Xi_{+} \left( T^* M \otimes \un(3)^{\perp} \right) = \Xi_{-}
\left( T^* M \otimes \un(3)^{\perp} \right) = T^* M \otimes T^* M
\wedge \omega.
\end{equation*}
Denote $(d\Psi_{\pm})_{4,5}$ the projections of $d\Psi_{+}$ and
$d\Psi_{-}$ onto the space $\mathcal W_{4,5}^a= T^* M \wedge
\Psi_+=T^* M \wedge \Psi_- \subseteq \Lambda^4T^*M $ respectively
defined by the alternating maps $\mathcal W_4^{\Xi} +\mathcal
W_5^{-} \longrightarrow \mathcal W_{4,5}^a$ and $\mathcal
W_4^{\Xi} +\mathcal W_5^{+} \longrightarrow \mathcal W_{4,5}^a$,
where $\mathcal W_4^{\Xi} = \Xi_+ (\mathcal W_4) = \Xi_- (\mathcal
W_4)$.

If we compute the $\mathcal W_4$-part $\left(\nabla \omega
\right)_4$  of $\nabla \omega$, the images $\Xi_{\pm} \left(\nabla
\omega \right)_4$   and then taking the skew-symmetric parts of
$\Xi_{\pm} \left(\nabla \omega \right)_4  \mp 3 J \eta \otimes
\Psi_{\mp}$, we will obtain the ${\mathcal W}_{4,5}$-parts of $d
\Psi_+$ and $d \Psi_-$, i.e.
\begin{equation} \label{dcuatrocinco}
(d \Psi_{\pm})_{4,5} = - \left( 3 \eta + \frac12 \theta^6 \right)
\wedge \Psi_{\pm}.
\end{equation}
With the help of   \eqref{estrella1}, \eqref{estrella2} and
\eqref{dcuatrocinco}, one gets that the one-form $\eta$  satisfies
the conditions
\begin{equation}\label{torsion3w5}
\ast \left( \ast  d \Psi_{\pm}  \wedge  \Psi_{\pm} \right) = 6
\eta +  \theta^6 = - J \ast \left( \ast  d \Psi_+ \wedge
\Psi_- \right) =    J \ast \left( \ast  d \Psi_- \wedge  \Psi_+
\right).
\end{equation}
So, we get the $\Lie{SU}(3)$-splitting \cite{CS}
\begin{gather*}
\eta + \xi^{\Un(3)} \in T^*M\otimes \lie{su}(3)^{\perp} =
{\mathcal W}^+_1\oplus{\mathcal W}^-_1 \oplus {\mathcal W}_2^+
\oplus{\mathcal W}_2^-\oplus {\mathcal W}_3 \oplus {\mathcal
W}_4\oplus\mathcal W_5 \subseteq T^*M \otimes \mbox{End} (TM).
\end{gather*}
Moreover, we have also
\begin{gather*}
\nabla \omega = - \xi^{\Un(3)} \omega \in T^*M\otimes
\lie{u}(3)^{\perp} = {\mathcal W}^+_1 \oplus{\mathcal W}^-_1
\oplus {\mathcal W}_2^+ \oplus{\mathcal W}_2^-\oplus {\mathcal
W}_3 \oplus {\mathcal
W}_4 \subseteq T^* M \otimes \Lambda^2 T^* M,\\
d \Psi_+, d \Psi_- \in {\mathcal W}_1^a \oplus {\mathcal W}_2^a
\oplus {\mathcal W}_{4,5}^a \subseteq \Lambda^4 T^* M,
\end{gather*}
where ${\mathcal W}_1^a = \mathbb R \omega \wedge \omega$, and
${\mathcal W}_2^a = \lie{su}(3) \wedge \omega$. Note that, using
the maps $\xi^{\Un(3)} \to - \xi^{\Un(3)} \omega=\nabla \omega$
and $\nabla \omega \to ( \mbox{Alt} \circ \Xi_{\pm} ) (\nabla
\omega)$, where $\mbox{Alt}$ denotes the alternation map, one has
the correspondences
\begin{gather*}
\left( \xi^{\Un(3)} \right)_{\mathcal W_j^+} \leftrightarrow
\left( \nabla \omega \right)_{\mathcal W_j^-}  \leftrightarrow
\left( d \Psi_+ \right)_{\mathcal W_j^a} = \mbox{Alt} \circ
\Xi_{+} \left( \nabla \omega
\right)_{\mathcal W_j^-}, \\
\left( \xi^{\Un(3)} \right)_{\mathcal W_j^-} \leftrightarrow
\left( \nabla \omega \right)_{\mathcal W_j^+}  \leftrightarrow
\left( d \Psi_- \right)_{\mathcal W_j^a} = \mbox{Alt} \circ
\Xi_{-} \left( \nabla \omega \right)_{\mathcal W_j^+}.
\end{gather*} 

We will also need an alternative approach to describe the summand
$\xi^{\Un(3)}$ of the intrinsic torsion  of an $\SU(3)$-structure.
We can write
\begin{equation} \label{covaseis}
\nabla \omega =-\xi^{\Un(3)} \omega= \sum_{j,k=1}^6 c_{jk} e_j
\otimes e_k \lrcorner \Psi_+.
\end{equation}
Consider the $\SU(3)$-map $r : T^* M \otimes \un(3)^{\perp} \to
\otimes^2 T^*M$ defined by
\begin{equation} \label{erresu}
 r(\beta)(x,y) = \frac{1}{2} \langle x
\lrcorner \beta , y \lrcorner \Psi_+   \rangle.
\end{equation}
It is straightforward to check that, for $\beta=\nabla\omega$ satisfying
 ~\eqref{covaseis}, $r(\nabla\omega) = \sum_{j,k=1}^6c_{ij} e_j
\otimes e_k$ and  the coderivative  $d^*\omega$ has the form
\begin{equation} \label{codeseis}
d^* \omega = \sum_{j=1}^{6}\;  \sum_{\{k,l |
\Psi_+(e_j,e_k,e_l)=1\}} \left(r(\nabla\omega)(e_k,e_l)-
r(\nabla\omega)(e_l,e_k)\right) e_j.
\end{equation}

A useful explicit description of the $\Lie{SU}(3)$-torsion $\eta +
\xi^{\Un(3)}$
 is presented in \cite{Cabrera:horient}. Since $\eta$ is given by \eqref{torsion3w5},
it remains to describe $\xi^{\Un(3)}$.  Write
$(d\Psi_{\pm})_{\xi^{\Un(3)}} = d \Psi_{\pm} + 3 \eta \wedge
\Psi_{\pm}$, and
 $(X\wedge Y)\lrcorner \left( d \Psi_{\pm} \right)_{\xi^{\Un(3)}} = \left( d\Psi_{\pm} \right)_{\xi^{\Un(3)}}(X,Y,\cdot,\cdot
)$. Then \cite{Cabrera:horient}
\begin{gather*} \label{exptorsion:xi}
 \hspace{-2.6cm} \xi^{\Un(3)}_{X} Y = - \frac{1}{2}\sum_{j,k=1}^6 r(\nabla\omega)(X, e_j) \Psi_{+}
(e_j,e_k,Y) Je_k,\\
 \label{erreextalg} 2 r(\nabla\omega)(X,Y) =
\left\langle X \lrcorner d \omega , Y \lrcorner \Psi_+
\right\rangle + \langle (JX \wedge Y) \lrcorner (d
\Psi_-)_{\xi^{\Un(3)}}- (X \wedge Y) \lrcorner (d
\Psi_+)_{\xi^{\Un(3)}} \, , \,\omega \rangle,
\end{gather*}
for all vectors $X,Y$.

The different classes of $SU(3)$-structures
can be characterized in terms  $d\omega$, $d \Psi_+$ and $d
\Psi_-$, as follows:
\begin{itemize}
\item{$\mathcal W_1 \oplus \mathcal W_5 =\mathcal
W_1^+\oplus\mathcal W_1^- \oplus \mathcal W_5 $:} The class of
nearly K\"ahler manifolds defined by $d\omega$ to be
(3,0)+(0,3)-form, i.e. $d \omega \in \mathbb R \Psi_+ \oplus
\mathbb R \Psi_-$, and $ d \Psi_{\pm} + 3 \eta \wedge \Psi_{\pm} \in
\mathbb R \omega \wedge \omega$.
 \item{$\mathcal W_2 \oplus \mathcal W_5 =\mathcal W_2^+\oplus\mathcal W_2^- \oplus \mathcal W_5$:}
The class of almost K\"ahler manifolds defined by $d\omega=0$.
\item{$\mathcal W_3 \oplus \mathcal W_5$:} The class of balanced
Hermitian manifolds determined by $d\Psi_{\pm} = \theta^6=0$.
 \item{$\mathcal W_4 \oplus \mathcal
W_5$:} The class of locally conformally K\"ahler spaces defined by
$2 d\omega= \theta^6\wedge\omega$.
  \item{$\mathcal W_5$:} The class of
K\"ahler spaces determined by the one-form $\eta$ given by
\eqref{torsion3w5}.
\end{itemize}
Note that if all components are zero, then we have a Ricci-flat
K\"ahler manifold. If the complex volume form is closed,
$d\Psi=0$, one gets the  observation due to Hitchin \cite{hitc}
that the almost complex structure is integrable.

A  new object is the class of \emph{half-flat} ( or $\mathcal
W_1^- \oplus \mathcal W_2^- \oplus \mathcal W_3$)
$\Lie{SU}(3)$-manifolds  which can be characterized by the
conditions
\begin{equation}\label{nw1}
d\Psi_+=\theta^6=0.
\end{equation}
The half-flat  $\SU(3)$-structures  can be lifted to a
$\Lie{G}_2$-holonomy metric on the product by the real line
solving the Hitchin flow equations \cite{Hit}. In fact, many new
$\Lie{G}_2$-holonomy metrics are obtained in this way
\cite{BGom,GLPS,CF,CGLP}.

\subsection{$\Spin(7)$-structures}

Now, let us consider $\mathbb R^8$ endowed with an orientation and
its standard inner product. Let $\{e,e_0,...,e_6\}$ be an oriented
orthonormal basis.  Consider the four-form $\Phi$ on $\mathbb R^8$
given by
\begin{eqnarray}\label{1}
\Phi &=& \sum_{i \in \mathbb{Z}_7}  e\wedge e_i \wedge e_{i+1}
\wedge  e_{i+3} - \sigma \sum_{i \in \mathbb{Z}_7} e_{i+2} \wedge
e_{i+4} \wedge e_{i+5} \wedge  e_{i+6},
\end{eqnarray}
where $\sigma$ is a fixed constant such that $\sigma = + 1$ or
$\sigma=-1$, and $+$ in the subindexes means the sum in
$\mathbb{Z}_7$.  We fix $e \wedge e_0 \wedge \dots \wedge e_6 = \frac{\sigma}{14} \Phi
\wedge \Phi$ as a volume form.

The subgroup of $\GL(8,\mathbb R)$ which fixes $\Phi$ is
isomorphic to the double covering $\Spin(7)$ of $\SO(7)$
\cite{HL}. Moreover, $\Spin(7)$ is a compact simply-connected Lie
group of dimension 21 \cite{Br}. The Lie algebra $\lie{spin}(7)$
of $\Spin(7)$ is isomorphic to the skew-symmtric two-forms $\psi$
satisfying the  linear equations
$$
\sigma \psi ( e_{i},e) + \psi ( e_{i+1}, e_{i+3}) + \psi (
e_{i+4}, e_{i+5}) + \psi ( e_{i+2}, e_{i+6}) = 0,
$$
for all $i \in \mathbb{Z}_7$. Shortly, $\lie{spin}(7)\cong  \{\psi
\in \Lambda^2 T*M |*_8( \psi\wedge\Phi)= \psi \}$.  The orthogonal
complement $\lie{spin}(7)^{\perp}$ of $\lie{spin}(7)$ in
$\Lambda^{2} \mathbb R^{8\ast} = \lie{so}(8)$ is the
seven-dimensional space generated by
\begin{equation}\label{new1}
\beta_i = \sigma e_{i} \wedge e + e_{i+1} \wedge e_{i+3} + e_{i+4}
\wedge e_{i+5} + e_{i+2} \wedge e_{i+6},
\end{equation}
where $i \in \mathbb{Z}_7$. Equivalently, $\lie{spin}(7)^{\perp}$
is described  as the space  consisting of those skew-symmetric
two-forms $\psi$ such that $*_8( \psi\wedge\Phi)= - 3 \psi$.

A $\Spin(7)$\emph{-structure} on an eight-manifold $M^8$ is by
definition a reduction of the structure group of the tangent
bundle to $\Spin(7)$; we shall also say that $M$ is a
\emph{$\Spin(7)$-manifold}. This can be geometrically described by
saying that there exists a nowhere vanishing global differential
four-form $\Phi$ on $M^8$ and a local frame $\{e,e_0,\dots,e_6\}$
such that the four-form $\Phi$  can be locally written as in
\eqref{1}. The four-form $\Phi$ is called the \emph{fundamental
form} of the $\Spin(7)$-manifold $M$ \cite{Bo} and the  local
  frame  $\{e, e_0, \dots, e_6 \}$ is called {\it a Cayley frame}.

The fundamental form of a $\Spin(7)$-manifold determines a
Riemannian metric $\langle \cdot , \cdot\rangle$ through
 $\langle x,y \rangle = - \frac17 \ast_8 \left( (x \lrcorner \Phi)
 \wedge \ast_8 \left( y \lrcorner \Phi \right) \right)$
 \cite{Gr}. Thus, $\langle \cdot , \cdot \rangle$ is referred as the metric induced by $\Phi$.
 Any  Cayley frame becomes an orthonormal frame with respect to such a metric.
We recall that the corresponding three-fold vector cross
  product $P$ is defined by
  $$
 \langle P(X_1,X_2,X_3),X_4 \rangle = \Phi (X_1,X_2,X_3,X_4),
  $$
  for smooth vector fields $X_i$ on $M^8$.

In general, not every eight-dimensional Riemannian spin manifold
$M^8$ admits a $\Spin(7)$-structure. We explain the precise
conditions given in \cite{LM}. Denote by $p_1(M)$, $p_2(M)$,
${\mathbb X}(M)$, ${\mathbb X}(S_{\pm})$ the first and the second
Pontrjagin classes, the Euler characteristic of $M$ and the Euler
characteristic of the positive and the negative spinor bundles,
respectively. It is well known \cite{LM} that a spin
eight-manifold admits a $\Spin(7)$-structure if and only if
${\mathbb X}(S_+)=0$ or ${\mathbb X}(S_-)=0$. The latter
conditions are equivalent to $ p_1^2(M)-4p_2(M)+ 8{\mathbb
X}(M)=0$, for an appropriate choice of the orientation.

Let us recall that a $\Spin(7)$-manifold $(M,\langle \cdot , \cdot
\rangle, \Phi)$ is said to be parallel (torsion-free),  if the
holonomy of the metric $Hol(\langle \cdot , \cdot \rangle)$ is a
subgroup of $\Spin(7)$. This is equivalent to saying that the
fundamental form $\Phi$ is parallel with respect to the
Levi-Civita connection $\nabla$ of the metric $\langle \cdot ,
\cdot \rangle$. Moreover, $Hol(\langle \cdot , \cdot
\rangle)\subseteq \Spin(7)$ if and only if $d\Phi=0$ \cite{F,Br}
(see also \cite{Sal}) and any parallel $\Spin(7)$-manifold is
Ricci-flat \cite{Bo}. The first known explicit example of complete
parallel $\Spin(7)$-manifold with $Hol(\langle \cdot , \cdot
\rangle)=\Spin(7)$ was constructed by Bryant and Salamon
\cite{BS,Gibb}. The first compact examples of parallel
$\Spin(7)$-manifolds with $Hol(\langle \cdot , \cdot
\rangle)=\Spin(7)$ were constructed by Joyce \cite{J1}.

There are four classes of $\Spin(7)$-manifolds according to
 Fern{\'a}ndez classification \cite{F}
obtained as irreducible $\Spin(7)$-representations   of the space
$\overline W \cong \mathbb R^{8*} \otimes \lie{spin}(7)^{\perp}$
of all possible covariant derivatives $\nabla\Phi$ of the
fundamental form with respect to the Levi-Civita connection. The
Lee form $\theta^8$ is defined by \cite{C1}
\begin{equation}\label{c2}
\theta^8 = -\frac{1}{7}*(*d\Phi\wedge\Phi)=\frac17*(\delta\Phi\wedge
\Phi).
\end{equation}
Fern{\'a}ndez classification can be described in terms of the Lee
form as follows : $\overline W_0 : d\Phi=0; \quad \overline W_1 :
\theta^8=0; \quad \overline W_2 : d\Phi = \theta^8 \wedge\Phi;
\quad \overline W :\overline W=\overline W_1\oplus \overline W_2$.

A $\Spin(7)$-structure of the class $\overline W_1$ (i.e.
 $\Spin(7)$-structure with
zero Lee form) is called
 \emph{a balanced $\Spin(7)$-structure}.
If the Lee form is closed, $d \theta^8=0$, then the
$\Spin(7)$-structure is locally conformal equivalent to a balanced
one \cite{I1}. It is shown in \cite{C1} that the Lee form of a
$\Spin(7)$-structure in the class $\overline W_2$ is closed.
Therefore,  such a manifold is locally conformal equivalent to a
parallel $\Spin(7)$-manifold. Compact spaces with closed but not
exact Lee form (i.e. the structure is not globally conformal
parallel) have very different topology than the parallel ones
\cite{I1}. Coeffective cohomology and coeffective numbers of
Riemannian manifolds with $\Spin(7)$-structure are studied in
\cite{Ug}.

\section{Intrinsic torsion of  $\Spin(7)$-structures}
In \cite{Br}, Bryant predicted  the existence of a formula
expressing the covariant derivative $\nabla\Phi$ of the
fundamental four-form in terms of its exterior derivative $d\Phi$
(see also \cite{Sal}). An explicit  expression of $\nabla \Phi$ in
terms of $d\Phi$ has been given in \cite{I1}. In this section we
use the alternative way of characterizing  the different types of
$\Spin(7)$-structure proposed in \cite{C3,Cabrera:hipspin}. This
help us to describe explicitly the intrinsic torsion of a given
$\Lie{Spin}(7)$-structure and to get a formula for $\nabla\Phi$ in
terms $d\Phi$. We note that the general properties of the
$\Lie{Spin}(7)$-intrinsic torsion are established in \cite{Cley}.

We consider the $\Lie{Spin}(7)$-isomorphism $\overline{r} \, : \,
\overline W \to \mathbb R^{8*} \otimes \lie{spin}(7)^{\perp}
\subset \mathbb R^{8*} \otimes \Lambda^2 \mathbb R^{8*}$  defined
by
$$
\overline{r}(\mathcal B)(x,y,z) = \frac1{8} \langle x \lrcorner
\mathcal B, y \wedge z \lrcorner \Phi - z \wedge y \lrcorner \Phi
\rangle, \quad x,y,z \in \mathbb R^{8}, \mathcal B \in \overline
W.$$
It is easy to see that $\overline{r}$ is a
$\Lie{Spin}(7)$-map. On the other hand, any $\mathcal B \in
\overline W$ can be written in the form \cite{C1}
\begin{equation}\label{ALDERsp}
\mathcal B =  \sigma \sum_{i \in {\mathbb Z}_{\;8},j \in {\mathbb
Z}_{\;7}} a_{ij} \,  e_{i} \otimes (e_j \wedge e \lrcorner \Phi -
e \wedge e_j \lrcorner \Phi ),
\end{equation}
where $\{ e=e_7, e_0, \dots , e_6 \}$ is a Cayley frame. Now one
can easily check that
\begin{equation}\label{ralfa1}
\overline{r}(\mathcal B) =  \sum_{i \in {\mathbb Z}_{\;8}, j \in
{\mathbb Z}_{\;7}} a_{ij} e_{i} \otimes \beta_{j},
\end{equation}
where the two-forms $\beta_j$ are determined in \eqref{new1}.
Therefore,  $\overline{r}$ is an isomorphism and the four classes
of $\Lie{Spin}(7)$-structures are expressed in terms of $\overline
r$ in \cite{C3}.

Further, we  describe the intrinsic $\Spin(7)$-torsion in terms of
$d \Phi$. Taking the skew-symmetric part of $\nabla\Phi$ given by
  \eqref{ALDERsp}, we obtain
\begin{eqnarray}\nonumber
d \Phi = & - \sum_{i \in {\mathbb Z}_{\;7}} \left( (
a_{i+2,i+2} + a_{i+4,i+4} + a_{i+5,i+5} + a_{i+6,i+6}) \, e \wedge
e_{i+2} \wedge e_{i+4} \wedge e_{i+5} \wedge e_{i+6} \right .
\\\nonumber
& +  ( \sigma \,  a_{7,i} + a_{i+4,i+5} + a_{i+1,i+3} +
a_{i+2,i+6})\, e \wedge e_{i}
\wedge e_{i+1} \wedge e_{i+2} \wedge e_{i+4} \\\nonumber
& +  (  \sigma \, a_{7,i} - a_{i+5,i+4} - a_{i+3,i+1} +
a_{i+2,i+6})\, e \wedge e_{i}
\wedge e_{i+2} \wedge e_{i+3} \wedge e_{i+5} \\\nonumber
& +  ( \sigma \, a_{7,i} + a_{i+4,i+5} - a_{i+3,i+1} -
a_{i+6,i+2})\, e \wedge e_{i}
\wedge e_{i+3} \wedge e_{i+4} \wedge e_{i+6} \\\label{diffund}
& +  (  \sigma \, a_{7,i} - a_{i+5,i+4} + a_{i+1,i+3} -
a_{i+6,i+2})\, e \wedge e_{i}
\wedge e_{i+5} \wedge e_{i+6} \wedge e_{i+1} \\\nonumber
& +\sigma \,   ( a_{i+4,i+5} - a_{i+5,i+4} + a_{i+1,i+3} -
a_{i+3,i+1})\, e_{i} \wedge e_{i+1}
\wedge e_{i+3} \wedge e_{i+4} \wedge e_{i+5} \\\nonumber
& + \sigma \, ( a_{i+4,i+5} - a_{i+5,i+4} + a_{i+2,i+6} -
a_{i+6,i+2})\, e_{i} \wedge e_{i+2}
\wedge e_{i+4} \wedge e_{i+5} \wedge e_{i+6} \\\nonumber
&\left. + \sigma \,      ( a_{i+1,i+3} - a_{i+3,i+1} + a_{i+2,i+6}
- a_{i+6,i+2})  \, e_{i} \wedge e_{i+6} \wedge e_{i+1} \wedge
e_{i+2} \wedge e_{i+3} \right). \nonumber
\end{eqnarray}
Consequently,  for the Lee form $\theta^8$, \eqref{diffund} and
\eqref{c2} yield
\begin{equation}
\label{pdelta}
\begin{array}{rcl}
\theta^8 & = &  - \frac47   \sum_{i \in {\mathbb Z}_{\;7}} (
a_{i+4,i+5} - a_{i+5,i+4} + a_{i+1,i+3} - a_{i+3,i+1} +
a_{i+2,i+6} - a_{i+6,i+2}
 + \sigma \, a_{7,i} ) \,  e_{i} \\[2mm]
 && + \,  \frac47 \sigma \sum_{i \in {\mathbb Z}_{\;7}}
a_{i,i} \, e.
\end{array}
\end{equation}
The equalities \eqref{diffund} and  \eqref{pdelta} imply
\begin{pro}
For a $\Lie{Spin}(7)$-structure, the condition $d\Phi=\theta^8
\wedge\Phi$  is equivalent to
\begin{equation} \label{rw2}  4 \overline{r} (\nabla\Phi) = \sum_{i \in \mathbb{Z}_8} e_i \otimes e_i
 \wedge\theta^8 + \sigma \theta^8 \lrcorner \Phi. \end{equation}
\end{pro}
Further, we have
\begin{thm} \label{spin7t}
 The minimal $\Spin(7)$-connection is given by
$\nabla^{\Spin(7)} = \nabla + \xi^{\Spin(7)}$, where the intrinsic
torsion $\xi^{\Spin(7)}$ is determined by
$$
 \langle \xi^{\Spin(7)}_X Y , Z \rangle = \frac{1}{4} \;
\overline{r}(\nabla\Phi)(X,Y,Z).
$$
Equivalently,
$$
\xi^{\Spin(7)}_X Y =  - \frac{\sigma}{24} \sum_{i , j \in {\mathbb
Z}_8} \overline{r}(\nabla\Phi) (X, e_i,e_j) P(e_i,e_j,Y),
$$
where $\{ e=e_7, e_0, \dots, e_6 \}$ is a Cayley frame.

The tensor $\overline{r}(\nabla\Phi)$ is
expressed in  terms of $d \Phi$  due to the next equality
\begin{equation} \label{rxyz}
4\overline{r}(\nabla\Phi)(X,Y,Z) =  2 \langle X \lrcorner d \Phi,  Y
\wedge Z\lrcorner \Phi -Z \wedge Y\lrcorner \Phi \rangle -  7(X
\wedge \theta^8)  (Y,Z).
\end{equation}
\end{thm}
\begin{proof} Let  $i \in \mathbb{Z}_8$ and $j \in \mathbb{Z}_7$. Then \eqref{ralfa1} and
\eqref{new1} give $4\overline{r}(\varphi)(e_i,e_j,e)= 4 \sigma
a_{ij}$. Now,  using the expressions (\ref{diffund}) and
(\ref{pdelta}) for $d\Phi$ and $\theta^8$, respectively, we check
that the right hand side of \eqref{rxyz} (denote it by $C$) gives
$C(e_i,e_j,e)= 4 \sigma
a_{ij}=4\overline{r}(\nabla\Phi)(e_i,e_j,e)$. Likewise, using
again  \eqref{diffund} and
  \eqref{pdelta}, one  checks that
$$
\sigma C(e_i,e_j,e) = C(e_i,e_{j+1} ,e_{j+3}) = C(e_i,e_{j+4},e_{j+5}) = C(e_i,e_{j+2} ,e_{j+6}).
$$
Therefore, $C \in T^* M \otimes \lie{spin}(7)^{\perp}$ and $4
\overline{r} (\Phi) = C$.
 In a similar way, one verifies that
$\xi^{\Lie{Spin}(7)} \in T^* M^8 \otimes \lie{spin}(7)^{\perp}$.
Finally, it is straightforward to check that
$\nabla^{\Lie{Spin}(7)} \Phi = 0$. Hence $\nabla^{\Lie{Spin}(7)}$
is a $\Lie{Spin}(7)$-connection.
\end{proof}
\begin{co}
The covariant derivative $\nabla \Phi$ of the fundamental form is expressed in terms of the exterior derivative
$d\Phi$  as follows $$\nabla \Phi = -
\xi^{\Lie{Spin}(7)} \Phi,$$
where $\xi^{\Lie{Spin}(7)}$ is determined in Theorem~\ref{spin7t}.
\end{co}

\section{$\SU(3)$-structures on six-dimensional submanifolds}

Let $f:M^6\longrightarrow (M^8,\Phi, \langle \cdot , \cdot
\rangle)$ be a smooth orientable six-manifold immersed in an
eight-dimensional $\Lie{Spin}(7)$-manifold with fundamental form
$\Phi$ and Riemannian metric $\langle \cdot , \cdot \rangle$.

Let $N_1,N_2$ be a local orthonormal frame of the normal bundle
$T^{\perp}M^6$. The $\Lie{Spin}(7)$-structure on $M^8$ induces an
almost Hermitian structure on $M^6$ defined \cite{Gr}
\begin{equation} \label{alhermstruc}
JX=P(N_1,N_2,X), \qquad X\in TM^6,
\end{equation}
where $P$ is the three-fold vector cross product on $M^8$ determined by the
$\Lie{Spin}(7)$-structure.

 It is well known that the almost complex structure $J$ is
independent on the particular oriented orthonormal frame and is
compatible with the induced Riemannian metric on $M^6$ \cite{Gr}.
Thus, we have a natural global almost Hermitian structure on
$M^6$, where  the K{\"a}hler form $\omega$ and  the Hodge star
operator $\ast_6$ are determined by $$\omega = \ast_6 \sigma f^*
\Phi, \qquad -4 \sigma Vol_6 = f^* \left( N_1 \lrcorner \Phi
\right) \wedge f^* \left( N_2 \lrcorner \Phi \right).$$  Also note
that $- 2 \sigma f^* \Phi = \omega \wedge \omega$.

As we have already pointed out, in general, there is not a  global
$\SU(3)$-structure induced from the $\Lie{Spin}(7)$-structure on
$M^8$. In fact, this assertion is based on the observation, due to
Bryant \cite{Br3}, saying that   the stabilizer of an oriented
two-plane in $\Lie{Spin}(7)$ is the group $\Lie{U}(3)$. In the
case $M^8=\mathbb R^8= \mathbb R^1\oplus Im{\mathbb O}$,  where
$Im{\mathbb O}$ is the space of imaginary octonions and
$M^6\subset Im{\mathbb O}$, there exists a global
$\Lie{SU}(3)$-structure due to the fact that the stabilizer in
$\Lie{Spin}(7)$ of two unitary vectors is the group $\Lie{SU}(3)$.
This phenomena  was discovered and studied by Calabi \cite{Cal}.
More general, any orientable hypersurface of a
$\Lie{G}_2$-manifold inherits a global $\SU(3)$-structure
\cite{Cal,Gray:someexamples,Cabrera:horient}.

We consider  local $\Lie{SU}(3)$-structures naturally induced from
the $\Lie{Spin}(7)$-structure on $M^8$. Namely, define the real
three-forms $\Psi_+, \Psi_-$ by the relations
\begin{gather}\label{defsu}
\Psi_+= \cos \gamma f^*(N_1\lrcorner\Phi) - \sin \gamma f^*(\sigma
N_2\lrcorner\Phi), \\\nonumber \Psi_-= \sin \gamma
 f^*(N_1 \lrcorner\Phi) + \cos \gamma f^*(\sigma N_2\lrcorner\Phi),
\end{gather}
where  $\gamma$ is a smooth function defined on $M^6$. The complex
three-form $\Psi$ with the real part $Re(\Psi)=\Psi_+$ and
imaginary part $Im(\Psi)=\Psi_-$ with respect to the induced
almost complex structure $J$ defined by \eqref{alhermstruc} is
clearly a local complex volume form compatible with the induced
$\Un(3)$-structure in the sense that it is a (3,0)-form with
respect to $J$, $J_{(1)} \Psi_{+} = \Psi_{-}$. Fixing $-\tfrac14
\Psi_+ \wedge \Psi_-$ as real volume form, the metric $\langle
\cdot , \cdot \rangle$ and the K\"ahler form $\omega$ are given by
$$
\langle x,y \rangle  = \frac{1}{2} \ast_6 \left( (x \lrcorner
\Psi_+ ) \wedge \ast_6 (y \lrcorner \Psi_+ )\right), \qquad \omega
(x,y) =   \frac{1}{2} \ast_6 \left( (x \lrcorner \Psi_- ) \wedge
\ast_6 (y \lrcorner \Psi_+ )\right),
$$
 respectively. The
three-forms $\Psi_+$ and $\Psi_-$ clearly depend on the local
orthonormal frame on the normal bundle. Therefore, they define a
local $\Lie{SU}(3)$-structure compatible with the global almost
Hermitian $\Lie{U}(3)$-structure $(\langle \cdot , \cdot
\rangle,J)$.

\begin{rem}\label{rem1} It is clear that all the local $\Lie{SU}(3)$-structures
generating the same metric are described by taking all oriented
orthonormal frames on the normal bundle and considering the
corresponding local $\Lie{SU}(3)$-structures defined above by
\eqref{defsu}. Also note that if we consider the local frame
$N_1'$, $N_2'$ on the normal bundle of $M^6$ given by $N_1' = \cos
\gamma N_1 - \sin \gamma N_2$ and $N_2'  = \sin \gamma N_1 + \cos
\gamma N_2$, then the complex volume form $\Psi$ defined in
\eqref{defsu} satisfy
  $\Psi_+ = f^*(N_1'\lrcorner\Phi)$ and $\Psi_- = f^*(\sigma N_2'\lrcorner\Phi)$. In this way we recover all local
 $\Lie{SU}(3)$-structures generating the same almost hermitian structure.
\end{rem}

The types of the induced global almost Hermitian
$\Un(3)$-structure depend on the second fundamental form of the
immersion and were described by Gray \cite{Gr} (see also
\cite{Br3}). We show below that the type of the induced local
$\SU(3)$-structures also depends  on the structure of the normal
bundle.

We  briefly recall some basic notions of the submanifold theory
(see e.g. \cite{Chen}).

Let us  fix an oriented orthonormal frame $N_1,N_2$ of the normal
bundle. Let $\nabla^8, \nabla^6$ be the Levi-Civita connection on
$M^8$, $M^6$, respectively. The Gauss equations read
\begin{equation}\label{gaus1}
\nabla^8_XY =\nabla^6_XY + \alpha(X,Y),\quad
\nabla^8_XN_j=-A_{N_j}X + D_XN_j, \quad j=1,2,\quad X,Y\in TM^6,
\end{equation}
where
\begin{equation}\label{gaus2}
\alpha(X,Y)=\alpha_1(X,Y)N_1 + \alpha_2(X,Y)N_2
\end{equation}
is the second fundamental form, $A_{N_j}$, $j=1,2$ is the shape
operator and $D$ is the normal connection. Since the normal
two-frame is orthonormal, we have
\begin{gather}\label{sh}
\langle A_{N_j}X,Y \rangle =\alpha_j(X,Y), \quad j=1,2, \quad
X,Y\in TM^6,\\\label{norm} D_XN_1=a(X)N_2, \qquad D_XN_2=-a(X)N_1,
\quad X\in TM^6,
\end{gather}
where $a(X)$ is a smooth function on $M^6$ depending on $X$.

When the shape operator vanishes, $M^6$ is said to be totally
geodesic.  The mean curvature $H$ is defined by $H=1/6 \, \Lie{tr}
\alpha=h_1N_1+h_2N_2$, where $6h_1=\Lie{tr} \alpha_1, \quad 6
h_2=\Lie{tr} \alpha_2$. The submanifold is said to be minimal, if
$H=0$, and totally umbilic, if $\alpha= \langle \cdot , \cdot
\rangle H$.
\subsection{Types of local $\SU(3)$-structures induced on six-dimensional
submanifolds} To investigate  special types of  local
$\Lie{SU}(3)$-structures, we find  relations between the local
intrinsic $\Lie{SU}(3)$-torsion of $M^6$ and  the global intrinsic
$\Lie{Spin}(7)$-torsion of the ambient manifold $M^8$. In the next
technical result,  we get relations involving the intrinsic
torsions, the shape operator and the structure of the normal
bundle of $M^6$.
\begin{pro} \label{rraa}
For the local $\Lie{SU}(3)$-structures on an oriented submanifold
$M^6$ of a $\Lie{Spin}(7)$-manifold $M^8$ inherited by the
$\Lie{Spin}(7)$-structure of $M^8$ and defined by \eqref{defsu},
we have the equalities
\begin{eqnarray} \label{rraa0}
r(\nabla^6\omega) & = & \cos \gamma (\sigma
\overline{r}(\nabla^8\Phi)(f_{\ast} \cdot,f_{\ast} J \cdot, N_1)
+ \sigma J_{(2)} \alpha_1 + \alpha_2) \\
&&  - \sin \gamma (  \sigma \overline{r}(\nabla^8\Phi)(f_{\ast}
\cdot ,f_{\ast} \cdot , N_1 ) - \sigma \alpha_1  + J_{(2)}
\alpha_2 ), \nonumber \\
 \theta^6& = &\frac74 f^*\theta^8  + \overline{r} (\nabla^8 \Phi
)(N_1, f_{\ast}  \cdot , N_1 )  - \sigma \overline{r} (\nabla^8
\Phi )(N_2, f_{\ast}  J \cdot, N_1 )  \label{rraa1} \\
&& + \sigma \overline{r} (\nabla^8 \Phi )(  f_{\ast} J \cdot ,
N_2,N_1),\nonumber \\
 \qquad \qquad \frac74 \; \sigma \theta^8 (N_1) & = &  \sigma
\overline{r} (\nabla^8 \Phi )(N_2,N_2,N_1) +  6 \sigma h_1    -
\sin \gamma \Lie{tr} r(\nabla^6\omega) +
2 \cos \gamma   \langle r(\omega), \omega \rangle  , \label{rraa2} \\
 - \frac74 \; \theta^8 (N_2) &  = &   \overline{r} (\nabla^8 \Phi
)(N_1,N_2,N_1) - 6  h_2+  \cos \gamma \Lie{tr} r(\nabla^6\omega) +
2 \sin \gamma \langle r(\omega)  , \omega \rangle , \label{rraa3}
\\
%
 3 \eta & = & - J \Lie{d} \gamma + \frac12 \ast_6 \left( \ast_6 f^*(
L_{N_1} \Phi) \wedge f^* (N_1 \lrcorner \Phi)
\right)- \overline{r}(\nabla^8\Phi) (N_1,  f_{\ast} \cdot , N_1 )   \label{rraa4}\\
&&
 - \sigma \overline{r}(\nabla^8\Phi) ( f_{\ast} J \cdot , N_2, N_1 ), \nonumber\\
 3 \eta & = &  - J
\Lie{d} \gamma + \frac12 \ast_6 \left( \ast_6 f^*( L_{N_2} \Phi)
\wedge f^* (N_2 \lrcorner \Phi) \right) + \sigma
\overline{r}(\nabla^8\Phi) (N_2, f_{\ast}J \cdot ,N_1)
  \label{rraa5}\\
&& - \sigma \overline{r}(\nabla^8\Phi) ( f_{\ast} J \cdot ,
N_2,N_1 ), \nonumber
\end{eqnarray}  \vspace{-12mm}

\begin{eqnarray}
\ast_6 \left( \ast_6 f^*( L_{N_1} \Phi) \wedge f^* (N_1 \lrcorner
\Phi) \right) & = & - \sigma J \ast_6 \left( \ast_6 f^*( L_{N_1}
\Phi) \wedge f^* (N_2 \lrcorner \Phi) \right), \label{rraa4p}\\
\ast_6 \left( \ast_6 f^*( L_{N_2} \Phi) \wedge f^* (N_2 \lrcorner
\Phi) \right) & = &  \sigma J \ast_6 \left( \ast_6 f^*( L_{N_2}
\Phi) \wedge f^* (N_1 \lrcorner \Phi) \right), \label{rraa5p}
\end{eqnarray}
where  $L$ denotes
Lie derivative.
\end{pro}
\begin{proof}
On any point of $M^6$, we  consider a Cayley frame  $\{
e=N_1 , e_0=N_2 , e_1 , \dots , e_6 \}$. Using   \eqref{ALDERsp}
and \eqref{ralfa1}, we obtain
$$
\sigma \overline{r} (\nabla^8\Phi) (e_i, e_1 , N_1 ) =  a_{i1} (\nabla^8_{e_i} \Phi)(N_1,N_2,e_4,e_6) =
\langle ( \nabla^8_{e_i}
P) (N_1,N_2, e_4), e_6 \rangle.
$$
$\,$ From these identities it is not hard to show
\begin{gather*}
\sigma \overline{r} (\nabla^8\Phi) (e_i,  e_1 , N_1 ) = -(\nabla^6_{e_i}
\omega)(e_4,e_6) + \sigma \alpha_1(e_i , e_1) + \alpha_2 (e_i ,
Je_1).
\end{gather*}
Since
$$
2 (\nabla^6_{e_i} \omega)(e_4,e_6) =  \langle  \nabla^6_{e_i}
\omega , Je_1 \lrcorner f^*(N_1\lrcorner \Phi)  \rangle = -
\langle \nabla^6_{e_i} \omega , e_1 \lrcorner f^*(\sigma
N_2\lrcorner \Phi) \rangle,
$$
we get
\begin{gather} \label{rraa1s}
\sigma \overline{r} (\nabla^8\Phi) (X, JY, N_1) =  \frac12 \langle
\nabla^6_X \omega ,  Y \lrcorner f^*(N_1\lrcorner \Phi)  \rangle +
\sigma \alpha_1 (X,JY) - \alpha_2 (X,Y), \\
\sigma \overline{r} (\nabla^8\Phi) (X,Y, N_1) =  \frac12 \langle
\nabla^6_X \omega , Y \lrcorner f^*(\sigma N_2\lrcorner \Phi)
\rangle + \sigma \alpha_1 (X,Y) + \alpha_2 (X,JY). \label{rraa2s}
\end{gather}
Now,   \eqref{rraa0} follows  from
\eqref{rraa1s} and \eqref{rraa2s}, using  \eqref{defsu} and \eqref{erresu}.

Next, we derive  \eqref{rraa1} from   \eqref{pdelta}, taking
\eqref{codeseis} and   \eqref{rraa0} for $\gamma=0$ into account.
Note that the Lee form $\theta^6$ is independent on the choice of
the complex volume form.

$\,$ From   \eqref{rraa0} we get
\begin{eqnarray}
\sigma \sum_{i=1}^6 \overline{r} (\nabla^8\Phi) (e_i,  e_i , N_1 )
- 6 \sigma h_1 = - \sin \gamma \Lie{tr} r(\nabla^6\omega) + 2 \cos
\gamma \langle r(\nabla^6\omega),
\omega \rangle, \label{traza1} \\
\sigma \sum_{i=1}^6 \overline{r} (\nabla^8\Phi) (e_i, J e_i , N_1
) + 6 h_2 = \cos \gamma \Lie{tr} r(\nabla^6\omega) + 2 \sin \gamma
\langle r(\nabla^6\omega), \omega \rangle. \label{traza2}
\end{eqnarray}
Now   \eqref{rraa2} and   \eqref{rraa3} follow from
\eqref{pdelta}, using   \eqref{traza1} and   \eqref{traza2}.

Take  $\gamma=0$. Then $\Psi_+=N_1 \lrcorner \Phi$ and $\Psi_- =
\sigma N_2 \lrcorner \Phi$. Apply  \eqref{torsion3w5} to get
\begin{gather} \label{aad1}
 \ast_6 \left( \ast_6 d f^*(N_1 \lrcorner \Phi) \wedge f^*
(N_1 \lrcorner \Phi) \right)= \ast_6 \left( \ast_6  d f^*(N_2
\lrcorner \Phi) \wedge f^* (N_2 \lrcorner \Phi) \right)=  \\
 = - \sigma  J\ast_6 \left( \ast_6  d f^*(N_1 \lrcorner \Phi) \wedge f^*
(N_2 \lrcorner \Phi) \right)= \sigma J  \ast_6 \left( \ast_6 d
f^*(N_2 \lrcorner \Phi) \wedge f^* (N_1 \lrcorner \Phi) \right).
\nonumber
\end{gather}
Use \eqref{estrella1}, \eqref{estrella2}, \eqref{aad1} and  \eqref{torsion3w5} for a generic $\gamma$ to
obtain
\begin{gather} \label{plusplus1}
 \ast_6 \left( \ast_6 d \Psi_+ \wedge \Psi_+ \right)= 6 \eta + \theta^6  = - 2 J d
 \gamma + \ast_6 \left( \ast_6  f^*d(N_1 \lrcorner \Phi) \wedge f^*
(N_1 \lrcorner \Phi) \right), \\
\ast_6 \left( \ast_6 d \Psi_+ \wedge \Psi_+ \right) =  6 \eta + \theta^6  = - 2 J d
 \gamma + \ast_6 \left( \ast_6  f^*d(N_2 \lrcorner \Phi) \wedge f^*
(N_2 \lrcorner \Phi) \right). \label{plusplus2}
\end{gather}
$\,$ From \eqref{diffund},   \eqref{pdelta} and \eqref{rraa1}, we
obtain
\begin{eqnarray*}
\ast_6 \left( \ast_6  f^*(N_1 \lrcorner d \Phi) \wedge f^* (N_1
\lrcorner \Phi) \right) & = & - \theta^6 - 2
\overline{r}(\Phi) (N_1, N_1, f_{\ast} \cdot ) + 2 \sigma
\overline{r}(\Phi) ( f_{\ast} J \cdot ,
N_2, N_1 ), \\
 \ast_6 \left( \ast_6  f^*(N_2 \lrcorner d \Phi) \wedge f^*
(N_2 \lrcorner \Phi) \right) & = &  - \theta^6 - 2 \sigma
\overline{r}(\Phi) (N_2, f_{\ast}J \cdot ,N_1) + 2 \sigma
\overline{r}(\Phi) ( f_{\ast} J \cdot , N_2,N_1 ),
\end{eqnarray*}
where we used the well known identity
\begin{equation}\label{new2}
d \left( N \lrcorner  \Phi \right) = L_N \Phi - N
\lrcorner d \Phi.
\end{equation}
Now,  \eqref{rraa4} and   \eqref{rraa5} follow from
\eqref{plusplus1} and   \eqref{plusplus2}. Finally, \eqref{rraa4p}
and    \eqref{rraa5p} are consequences of \eqref{aad1},
\eqref{diffund} and    \eqref{pdelta}, taking the identity
\eqref{new2} into account.
\end{proof}

Proposition~\ref{rraa} gives us chance to find relations between
the $\Lie{Spin}(7)$-structure on the ambient eight-dimensional
manifold and the local
 $\SU(3)$-structure inherited on the six-dimensional submanifold involving  the second
fundamental form.
\begin{thm} \label{wzerogeneral}
Let $M^8$ be an eight-dimensional Riemannian manifold with a parallel
 $\Lie{Spin}(7)$-structure. Let $M^6$ be an oriented six-dimensional submanifold
 of $M^8$ with the local
$\SU(3)$-structure defined by   \eqref{alhermstruc},
\eqref{defsu}. Then $M^6$ is  of type  $ {\mathcal  W}^+_1 \oplus
{\mathcal  W}_1^- \oplus {\mathcal  W}^+_2 \oplus {\mathcal W}_2^-
\oplus {\mathcal
 W}_3 \oplus {\mathcal W}_5$ and the following identities hold
\begin{gather} \label{wzero1}
 \ast_6 \left( \ast_6 f^*( L_{N_1} \Phi) \wedge f^* (N_1 \lrcorner \Phi)
 \right) =
\ast_6 \left( \ast_6 f^*( L_{N_2} \Phi) \wedge f^* (N_2 \lrcorner \Phi) \right)= \\
\hspace{1.63cm} = - \sigma J \ast_6 \left( \ast_6 f^*( L_{N_1}
\Phi) \wedge f^* (N_2 \lrcorner \Phi) \right)  =  \sigma J \ast_6
\left( \ast_6 f^*( L_{N_2} \Phi) \wedge f^* (N_1 \lrcorner \Phi)
 \right). \nonumber
\end{gather}
The precise conditions which characterized the types of local $\SU(3)$-structures on $M^6$ are
displayed in Table~\ref{tab:wzerogeneral}.

In particular:
\begin{itemize}
\item[a)] $M^6$ is a minimal submanifold if and only if the global
$\Un(3)$-structure belongs to the class $\mathcal
W_2\oplus\mathcal W_3$ in the Gray-Hervella classification.
\item[b)] The global $\Un(3)$-structure on $M^6$ is nearly
K\"ahler (type $\mathcal W_1$)
 if and only if the submanifold is totally umbilical.
  \item[c)] The global
$\Un(3)$-structure on $M^6$ is K\"ahler
 if and only if the submanifold is totally geodesic.
\end{itemize}
\end{thm}
\begin{proof}
The identities \eqref{wzero1} are direct consequences  of
\eqref{aad1}, \eqref{new2} and the condition $d\Phi=0$. Observe
that the latter implies $\overline{r}(\nabla^8\Phi)=0$. Now,
Table~\ref{tab:wzerogeneral} and the remaining part of Theorem
\ref{wzerogeneral} are consequences of the equations given in
 Proposition~\ref{rraa}.
\end{proof}

\begin{rem}
Note that  Theorem~\ref{wzerogeneral} includes the results
obtained by Gray in \cite{Gr}.
\end{rem}

\begin{thm} \label{wonegeneral}
Let $M^8$ be an eight-dimensional Riemannian manifold with a
 $\Lie{Spin}(7)$-structure having zero Lee form, $\theta^8=0$.
 Let $M^6$ be an oriented six-dimensional submanifold
 of $M^8$ with the local
$\SU(3)$-structures defined by   \eqref{alhermstruc},
\eqref{defsu}. Then:
\begin{itemize}
\item[a)] The precise conditions characterizing the types of the
local $\SU(3)$-structure are given in Table~\ref{tab:wonegeneral}.
\item[b)] The following identities hold
$$
\quad \; \begin{array}{l} \theta^6 =- \overline{r} (\nabla^8 \Phi
)(N_1,N_1, f_{\ast} \cdot ) - \sigma \overline{r} (\nabla^8 \Phi
)(N_2, f_{\ast}  J \cdot, N_1 )
 +\sigma \overline{r} (\nabla^8 \Phi )(f_{\ast}  J \cdot, N_2, N_1
 ),\\[1mm]
\Lie{tr} r(\nabla^6\omega) = - 2 \sin \gamma \left( \sigma h_1 -
\sigma \overline{r} (\nabla^8 \Phi )(N_2,N_2,N_1) \right) + 2 \cos
\gamma  \left(
 h_2 -  \overline{r} (\nabla^8 \Phi )(N_1,N_2,N_1) \right),\\[1mm]
  \langle r(\nabla^6\omega), \omega \rangle =    \cos \gamma \left( \sigma h_1 - \sigma
\overline{r} (\nabla^8 \Phi )(N_2,N_2,N_1) \right) +  \sin \gamma \left(
 h_2 -  \overline{r} (\nabla^8 \Phi )(N_1,N_2,N_1) \right).
\end{array}
$$
\end{itemize}
\end{thm}
\begin{proof} Using $\theta^8=0$, the equalities in Proposition~\ref{rraa} imply the assertion.
\end{proof}

\begin{thm} \label{wtwogeneral}
Let $M^8$ be an eight-dimensional Riemannian manifold with a
locally conformal parallel $\Lie{Spin}(7)$-structure, i.e.
$d\Phi=\theta^8\wedge\Phi$.
 Let $M^6$ be an oriented six-dimensional submanifold
 of $M^8$ with the local
$\SU(3)$-structures defined by   \eqref{alhermstruc},
\eqref{defsu}. Then:
\begin{itemize}
\item[a)]  The following identities hold
\begin{equation} \label{wtwo}
\begin{array}{l}
  \ast_6 \left( \ast_6 f^*( L_{N_1} \Phi) \wedge f^* (N_1
\lrcorner \Phi) \right) = \ast_6 \left( \ast_6 f^*( L_{N_2}
\Phi) \wedge f^* (N_2 \lrcorner \Phi) \right)= \\[1mm]
 \hspace{1.6cm} =
- \sigma J \ast_6 \left( \ast_6 f^*( L_{N_1} \Phi) \wedge f^* (N_2
\lrcorner \Phi) \right) =  \sigma J \ast_6 \left( \ast_6 f^*(
L_{N_2} \Phi) \wedge f^* (N_1 \lrcorner \Phi)
 \right),\\[2mm]
4 r(\nabla^6\omega) = \cos \gamma (  \sigma \theta^8(N_1) \omega
+  \Phi ( \theta^8 , N_1 , f^* \cdot , f^* J \cdot)
+ 4 \sigma J_{(2)} \alpha_1 + 4 \alpha_2) \\[1mm]
 \qquad \; \;\; \; \; \; \; \; \; \;  - \sin \gamma ( \sigma \theta^8(N_1) \langle \cdot ,
 \cdot \rangle +  \Phi (\theta^8 , N_1 , f^* \cdot , f^*  \cdot)
- 4\sigma \alpha_1 + 4 J_{(2)} \alpha_2 ), 
\\ [2mm]
\theta^6 =  f^*\theta^8, 
\\[2mm]
6 \eta = - 2 J  d \gamma + \ast_6 \left( \ast_6 f^*( L_{N_1} \Phi)
\wedge f^* (N_1 \lrcorner \Phi) \right)   + J d^* \omega,\\[2mm]
 \displaystyle  \frac{2}{3} \Lie{tr} r(\nabla^6\omega) = \sin
\gamma \left( 4 \sigma h_1 - \sigma \theta^8(N_1) \right) + \cos
\gamma \left(
 4 h_2 - \theta^8(N_2) \right),\\[3mm]
\displaystyle  \frac{4}{3} \langle r(\nabla^6\omega), \omega
\rangle = - \cos \gamma \left( 4 \sigma h_1 - \sigma \theta^8(N_1)
\right) + \sin \gamma \left(
 4 h_2 - \theta^8(N_2) \right).
 \end{array}
 \end{equation}
\item[b)] The precise conditions characterizing the types of the
local $\SU(3)$-structure are given in Table~\ref{tab:wtwogeneral}.
In particular:
\begin{itemize}
 \item[i)] The global $\Un(3)$-structure is locally conformal equivalent
to a nearly K\"ahler structure if and only if $M^6$ is totally
umbilic submanifold. If moreover $\theta^8$ is normal to $M^6$,
then the structure is nearly K\"ahler.
 \item[ii)]The global
$\Un(3)$-structure is locally conformal  K\"ahler  if and only if
$M^6$ is totally umbilic submanifold such that
$h_1=\frac14\theta^8(N_1), \quad h_2=\frac14\theta^8(N_2)$. If
moreover $\theta^8$ is normal to $M^6$, then it is a K\"ahler
structure.
\end{itemize}
\end{itemize}
\end{thm}
\begin{proof}
Since the  $\Lie{Spin}(7)$-structure is locally conformal
parallel, the  equality \eqref{rw2} is valid and $d\theta^8=0$.
Therefore, the equalities in a) as well as the conditions in
Table~\ref{tab:wtwogeneral} are direct consequences of \eqref{rw2}
and  Proposition~\ref{rraa}. The totally umbilical conditions are
derived in the same way as in the proof of the
Theorem~\ref{wzerogeneral}. Now i) follows from the recent result
\cite{Butru} which states that any six-dimensional almost
Hermitian manifold in the class $\mathcal W_1 \oplus \mathcal W_4$
is locally conformal to a nearly K\"ahler space. Finally, if
$\theta^8$ is normal to $M^6$, then \eqref{wtwo} shows that the
Lee form on $M^6$ vanishes.
\end{proof}
\begin{co}
 \label{wtwogeneraltan} Let $M^8$ be an eight-dimensional Riemannian
manifold with a locally conformal parallel
$\Lie{Spin}(7)$-structure, i.e. $d\Phi=\theta^8\wedge\Phi$.
 Let $M^6$ be an oriented six-dimensional submanifold
 of $M^8$ with the local
$\SU(3)$-structures defined by   \eqref{alhermstruc},
\eqref{defsu}. If the Lee form $\theta^8$ is tangent to $M^6$,
then the precise conditions characterizing the types of the local
$\SU(3)$-structure are given in Table~\ref{tab:wtwogeneraltan}. In
particular:
\begin{itemize}
 \item[i)]The global $\Un(3)$-structure is   K\"ahler if and only if $M^6$
is totally geodesic and the restriction of the Lee form to $M^6$
vanishes, i.e. $f^* \theta^8=0$.
 \item[ii)]The global
$\Un(3)$-structure is locally conformal  K\"ahler  if and only if
$M^6$ is totally geodesic.
 \item[iii)]The global $\Un(3)$-structure is of
type $\mathcal W_2 \oplus  \mathcal W_3 \oplus \mathcal W_4$  if
and only if $M^6$ is minimal.
\end{itemize}
\end{co}

\section{Holomorphic complex volume form}\label{I1}

We investigate the case when the induced local complex volume form
is closed, which implies, in particular, that the almost complex
structure is integrable \cite{hitc}.

We begin with
\begin{pro}\label{clos}
Let $(M^8,\Phi,g)$ be a $\Lie{Spin}(7)$-manifold.
 Let $M^6$ be an oriented six-dimen\-sio\-nal submanifold and
let $N_1$,  $N_2$ be any orthonormal frame of the normal bundle.
Then the complex volume form
 $$\Psi = \Psi_+ + i \Psi_-, \quad \Psi_+ = f^*(N_1\lrcorner\Phi), \quad \Psi_- = f^*(\sigma
N_2\lrcorner\Phi)$$ is closed, $d\Psi=0$, if an only if the next
two conditions hold simultaneously
\begin{equation}\label{clos1}
L_{N_1}\Phi|_{M^6} =(N_1\lrcorner d\Phi)|_{M^6}, \qquad
L_{N_2}\Phi =(N_2\lrcorner d\Phi)|_{M^6}.
\end{equation}
In particular, the almost complex structure is integrable.

If the $\Spin(7)$-structure is parallel, $d\Phi=0$, then the
complex volume form is closed exactly when
\begin{equation}\label{clos2}
L_{N_1}\Phi|_{M^6}=L_{N_2}\Phi|_{M^6}=0.
\end {equation}
In particular, if the normal bundle is parallel along the
submanifold, then there exists a local closed complex volume form
compatible with the induced global almost Hermitian
$\Un(3)$-structure.
\end{pro}
\begin{proof}
Take the exterior derivative in \eqref{defsu} and use \eqref{new2}
to get \eqref{clos1} and, consequently, \eqref{clos2}. The
integrability of the almost complex structure in the case of
closed complex volume form follows from the result of Hitchin
\cite{hitc}.

The Lie derivative is expressed in terms of the Levi-Civita
connection as follows
\begin{gather}\label{lie}
(L_N\Phi)(X,Y,Z,V)=(\nabla^8_N\Phi)(X,Y,Z,V)+\\\nonumber
\Phi(\nabla^8_XN,Y,Z,V) +\Phi(X,\nabla^8_YN,Z,V)
+\Phi(X,Y,\nabla^8_ZN,V) +\Phi(X,Y,Z,\nabla^8_VN).
\end{gather}
Since the normal bundle is parallel along $M^6$, we may choose a
parallel oriented normal two-frame. Take the corresponding complex
volume form, we see that it is closed due to \eqref{clos2} and
\eqref{lie}.
\end{proof}
As a consequence of the proof of Proposition~\ref{clos}, we get a
result which second part is essentially established in
\cite{Cabrera:horient}.
\begin{thm}\label{halfflat6}
There exist a local half-flat  $\SU(3)$-structure induced on a
six-dimensional submanifold of a parallel $\Spin(7)$-manifold if
and only if there exists a normal vector field which preserves the
parallel $\Spin(7)$-form restricted to the submanifold.

In particular, any orientable hypersurface $M^6\subset \mathbb
R^7=Im{\mathbb O}\subset{\mathbb O}$ carries a global  half-flat
$\SU(3)$-structure.
\end{thm}
\begin{proof}
Since $M^6\subset \mathbb R^7=Im{\mathbb O}\subset{\mathbb O}$, we
may take $\cos\gamma=1$ and 
$\nabla^8N_1=0$. Therefore, $d\Psi_+=0$  according to the proof of
Proposition~\ref{clos}. Hence, \eqref{nw1} are satisfied since
$\theta^6=0$.
\end{proof}

\subsection{Application to Calabi and Bryant examples}
Now we restrict our attention to the case $M^8={\mathbb O}$
studied in detail by Bryant in \cite{Br3}. In this case (even more
general, when the $\Lie{Spin}(7)$-structure of the ambient
manifold is parallel),  some of the $\Lie{U}(3)$-components of the
induced almost Hermitian structure are described by Gray \cite{Gr}
(see also \cite{Br3}). He showed that  the Lee form $\theta^6$ is
always zero and the submanifold $M^6$ is necessarily minimal.
Therefore, if the almost complex structure is integrable, then it
is balanced (type $\mathcal{W}_3$). Submanifolds with balanced
almost Hermitian structure are investigated by Bryant in
\cite{Br3}. He shows that if $M^6\subset{\mathbb O}$ inherits
complex and non-K{\"a}hler structure, then $M^6$ is foliated by
four-planes in ${\mathbb O}$ in a unique way, he calls this
foliation asymptotic ruling. He also obtains that if the
asymptotic ruling is parallel, then $M^6$ is a product of a fixed
associative four-plane $Q^4$ in ${\mathbb O}$ with a minimal
surface in the orthogonal four-plane. Moreover, Bryant found that
the Calabi examples, described in \cite{Cal}, are exactly those
complex $M^6$ with parallel asymptotic ruling which lie in
$Im{\mathbb O}\subset {\mathbb O}$, i.e. the minimal surface lies
in an associative three-plane in $Im {\mathbb O}$.

We investigate below when the local $\Lie{SU}(3)$-structures is
holomorphic in the case of parallel asymptotic ruling.

To be more precise, we explain the Bryant construction. Let
$\mathbb R^8={\mathbb O}=\mathbb R^4 \oplus Q^4$ be an orthogonal
sum of Cayley planes and let $S\subset \mathbb R^4$ be a surface.
Then $S\times Q^4\subset {\mathbb O}$ inherits a complex structure
if and only if $S$ is minimal in $\mathbb R^4$ and non-K{\"a}hler
provided $S$ is not a complex curve in $\mathbb R^4$ for some of
$\mathbb R^{4^,}$s complex structures \cite{Br3}. We have
\begin{thm}\label{cal}
Let $S \subset \mathbb R^4$ be a minimal surface in $\mathbb R^4$
such that $M^6=S \times Q^4\subset{\mathbb O}$ is a non-K{\"a}hler
complex manifold with respect to the $\Un(3)$-structure  induced
from ${\mathbb O}$. There exists a local holomorphic
$\SU(3)$-structure compatible with the $\Un(3)$-structure if and
only if $S$ is a minimal surface in a three-plane $\mathbb R^3$.
In this case the $SU(3)$-structure is globally defined and the
holomorphic volume form is parallel with respect to the Bismut
connection. In particular, the $\SU(3)$-structure described by
Calabi is holomorphic CYT structure.
\end{thm}
\begin{proof}
We need information for the Lie derivative of the fundamental
four-form in the normal direction due to Proposition~\ref{clos}.

Let us  fix an oriented orthonormal frame $N_1,N_2$ in the normal
bundle $T^{\perp}S\subset \mathbb R^4$ in $\mathbb R^4$ and a
local frame $X_3,X_4$ of the tangent bundle $TS$. We denote
$e_5,e_6,e_7,e_8$ the vectors in $Q^4$. We may write \eqref{sh}
and \eqref{norm} in the form
\begin{gather}\label{sh1}
A_{N_1}X_3 = \alpha_1(X_3,X_3)X_3 + \alpha_1(X_3,X_4)X_4, \qquad
A_{N_2}X_3 =\alpha_2(X_3,X_3)X_3 +
\alpha_2(X_3,X_4)X_4,\\
 \nonumber A_{N_1}X_4 = \alpha_1(X_4,X_3)X_3
+ \alpha_1(X_4,X_4)X_4, \qquad A_{N_2}X_4 =\alpha_2(X_4,X_3)X_3 +
\alpha_2(X_4,X_4)X_4.
\end{gather}
\begin{gather}\label{norm2} D_{X_3}N_1=a(X_3)N_2, \qquad
D_{X_3}N_2=-a(X_3)N_1,\\\nonumber D_{X_4}N_1=a(X_4)N_2, \qquad
D_{X_4}N_2=-a(x_4 )N_1.
\end{gather}
The minimality condition implies the equalities
\begin{equation}\label{min}
\alpha_1(X_3,X_3) + \alpha_1(X_4,X_4) =0, \qquad \alpha_2(X_3,X_3)
+ \alpha_2(X_4,X_4) =0.
\end{equation}
Using  \eqref{sh1}, \eqref{norm2}, we obtain from \eqref{lie} that
$(L_{N_j}\Phi)(X_k,e_l,e_m,e_p)=0$, for $j=1,2$, $k=3,4$ and
$l,m,p=5,6,7,8$,  since $Q^4$ is a Cayley four-plane. It remains
to investigate the case when two of the four vectors are tangent
to $S$. We need in addition to take into account  the minimality
condition \eqref{min}. We obtain
\begin{gather}\label{min1}
(L_{N_1}\Phi)(X_3,X_4,e_l,e_m) = a(X_3)\Phi(N_2,X_4,e_l,e_m) -
a(X_4)\Phi(N_2,X_3,e_l,e_m),\\\nonumber
(L_{N_2}\Phi)(X_3,X_4,e_l,e_m) = - a(X_3)\Phi(N_1,X_4,e_l,e_m) +
a(X_4)\Phi(N_1,X_3,e_l,e_m).
\end{gather}
Taking into account that $Q^4$ is a Cayley submanifold, we get
from \eqref{min1} that $L_{N_1}\Phi|_{M^6}=L_{N_2}\Phi|_{M^6}=0$
if and only if $a(X_3)=a(X_4)=0$, i.e. the normal connection is
flat. Now, Proposition~\ref{clos} and Remark~\ref{rem1} yield that
there is a local holomorphic complex volume form compatible with
the induced metric exactly when the minimal surface $S$ has flat
normal bundle. It is known  that a minimal submanifold of an
Euclidean space has flat normal connection if and only if it lies
in a three-dimensional plane $\mathbb R^3$ (see e.g. \cite{Chen}).
In this case, $\theta^6=d\Psi_+=d\Psi_-=0$. Apply Theorem~4.1 of
\cite{II} to conclude that the corresponding Bismut connection
preserves the complex volume form $\Psi$, i.e. it has holonomy
contained in $\Lie{SU}(3)$. Therefore,  the structure is
Calabi-Yau with torsion which completes the proof.
\end{proof}
Applying \cite[Theorem 4.1 ]{II}, we obtain in view of
Theorem~\ref{cal}
\begin{thm}\label{cal1}
Let $S\subset \mathbb R^4$ be a minimal surface in $\mathbb R^4$
such that $M^6=S \times Q^4\subset{\mathbb O}$ is a non-K{\"a}hler
complex manifold with respect to the $\Un(3)$-structure  induced
from ${\mathbb O}$. Then the Bismut connection of this
$\Un(3)$-structure has holonomy contained in $\SU(3)$ if and only
if $S$ is a minimal surface in a three-plane $\mathbb R^3$.

In particular, the holonomy of the Bismut connection of the
$\SU(3)$-structure described by Calabi is contained in $\SU(3)$.
Consequently, the compact complex non-K{\"a}hler six-manifolds
with holomorphically trivial canonical bundle constructed by
Calabi are balanced CYT-manifolds with respect to the  Calabi's
$\SU(3)$-structure.
\end{thm}

\section{Examples}

\begin{example}\label{exam1}
$ S^3 \times S^3$. Let us consider $\mathbb R^8$ with its standard
parallel $\Spin(7)$-structure. Thus, if $(x,x_0, \ldots, x_6)$ are
the global coordinates of $\mathbb R^8$, the $\Spin(7)$-structure
on $\mathbb R^8$ is the one such that $ \left\{
\frac{\partial}{\partial x}, \frac{\partial}{\partial x_0},
\ldots, \frac{\partial}{\partial x_6} \right\} $ is a Cayley
frame. For sake of simplicity, we will denote $e=
\frac{\partial}{\partial x}$ and $e_i=\frac{\partial}{\partial
x_i}$, for $i \in \mathbb Z_7$.

Let $S^3_{\;1} \times S^3_{\;2}$ be the six-submanifold of
$\mathbb R^8$ consisting of the product of two three-dimensional
spheres $ S^3_{\;1} \subseteq \left( \mathbb R^4 \right)_1 = span
\left\{ e,  e_0 ,  e_1 , e_3\right\}$ and $S^3_{\;2}  \subseteq
\left( \mathbb R^4 \right)_2 = span \left\{ e_2 , e_4  , e_5 , e_6
\right\}$.
 Fixing the oriented normal frame
$N_1= x e + x_0 e_0 + x_1 e_1 + x_3 e_3$, $N_2= x_2 e_2 + x_4 e_4
+ x_5 e_5 + x_6 e_6$, we consider the $\SU(3)$-structure on
$S^3_{\;1} \times S^3_{\;2}$ defined by  \eqref{alhermstruc} and
\eqref{defsu}. This $\SU(3)$-structure is globally defined on
$S^3_{\;1} \times S^3_{\;2}$,  since the stabilizer of two
orthonormal vectors in $\Spin(7)$ is the group $\SU(3)$ and is
compatible with the standard product metric on $S^3\times S^3$.

The tangent bundle of $S^3_{\;1} \times S^3_{\;2}$ is decomposed
into $ T \left( S^3_{\;1} \times S^3_{\; 2}\right) = T S^3_{\;1}
\oplus T S^3_{\;2}$ and,  for all $X \in T \left( S^3_{\;1} \times
S^3_{\;2} \right)$, we have the corresponding decomposition
$X=X_1+X_2$. The observation
 \begin{gather*}
 P (N_1,N_2,e),  P (N_1,N_2,e_0) ,  P (N_1,N_2,e_1) ,  P
 (N_1,N_2,e_3)\in T \left(\mathbb  R^4 \right)_2, \\
 P (N_1,N_2,e_2),  P (N_1,N_2,e_4) ,  P (N_1,N_2,e_5) ,  P
 (N_1,N_2,e_6)\in T \left(\mathbb R^4 \right)_1
\end{gather*}
 yields $J \left( T_p
S^3_{\;1} \right)  = T_p  S^3_{\;2}$ and $J \left( T_p S^3_{\;2}
\right) = T_p S^3_{\;1} $, for any point $p\in S^3_{\;1} \times
S^3_{\; 2} $.

The second fundamental form is given by $ \alpha_1 ( X, Y ) = -
\langle X_1 , Y_1 \rangle$,  $\alpha_2 ( X, Y ) = - \langle X_2 ,
Y_2 \rangle$. Consequently,   $ (1+J) \alpha_1 = \alpha_1 +
\alpha_2 =  2 h_1 \langle \cdot , \cdot \rangle$, $\; (1+J)
\alpha_2 = \alpha_2 + \alpha_1 = 2 h_2 \langle \cdot , \cdot
\rangle$.
 Using the results in Theorem \ref{wzerogeneral} and
Table \ref{tab:wzerogeneral}, we conclude that the
$\Lie{SU}(3)$-structure on $ S^3_{\;1} \times S^3_{\;2}$ is of
type $\mathcal W^+_1 \oplus \mathcal W^-_1 \oplus \mathcal W_3
\oplus \mathcal W_5$.

We  describe the $\mathcal W_5$-part $\eta$ of the intrinsic
$\SU(3)$-torsion.
 The Lie derivative $L_{N_1} \Phi$
restricted to $ S^3_{\;1} \times S^3_{\;2}$ is given by
\begin{gather*}
f^* \left( L_{N_1} \Phi \right)  = - Alt \langle \nabla^8_{ \cdot}
N_1 , P(\cdot, \cdot,  \cdot) \rangle =  Alt \left( \alpha_1 (
\cdot , P(\cdot,\cdot, \cdot) \right)  = -  2 \sigma
\Phi_{|T\left(S^3_{\;1}\times S^3_{\;2} \right)} =  \omega \wedge
\omega.
\end{gather*}
This can be checked using a Cayley frame $\left\{ N_1, N_2, u_1,
\dots , u_6 \right\}$, where $u_1,u_2,u_4 \in T S^2_{\;1}$ and
$u_3,u_5,u_6 \in T S^2_{\;2}$. Such a Cayley frame do exist
because the almost complex structure $J$ maps the tangent space of
one $S^3$ to the tangent space of the another $S^3$. Note also
that $\ast_6 f^* \left( L_{N_1} \Phi \right) = - 2 \omega$.
Since
$f^* ( N_1 \lrcorner \Phi)$ is a linear combination of $\Psi_+$
and $\Psi_-$, we have
$$
\ast_6 f^* \left( L_{N_1} \Phi \right) \wedge  f^* ( N_1 \lrcorner
\Phi) = - 2 \omega \wedge f^* ( N_1 \lrcorner \Phi) =0.
$$
Now, using Equation \eqref{rraa4}, we get $3 \eta = - J d \gamma$.
Hence, the $\mathcal W_5$-part, $\eta$, of the intrinsic
$\Lie{SU}(3)$-torsion vanishes exactly when $\gamma$ is a
constant.

We compute  the exterior derivatives $d \omega$, $d \Psi_+$ and $d
\Psi_-$. Consider three orthonormal vector fields $v_1,v_2,v_3$ in
$T\left(S^3\right)_1$ such that $\Phi (N_1 , v_1 , v_2, v_3)=1$
(or $v_3 = P(N_1,v_1,v_2)$). We know that $Jv_1,Jv_2,Jv_3 \in T
\left(S^3\right)_2$. Taking into account the expression for $\Psi$
given by \eqref{defsu}, we obtain $\Psi (v_1,v_2,v_3) =
e^{i\gamma}$. Therefore, $v_1,v_2,v_3,Jv_1,Jv_2,Jv_3$ is an
adapted basis for the $\Lie{U}(3)$-structure but not for the
$\Lie{SU}(3)$-structure considered. However, if we write,  for
$i=1,2,3$,
\begin{equation}\label{mnnov1}
u_i = e^{-i\frac{\gamma}{3}} v_i=  \cos \frac{\gamma}{3} \; v_i -
\sin \frac{\gamma}{3} \; Jv_i,
\end{equation}
then we have $\Psi (u_1,u_2,u_3)=1$ and hence
 $u_1,u_2,u_3,Ju_1,Ju_2,Ju_3$ is a local frame
adapted to the $\Lie{SU}(3)$-structure. For the second fundamental
form we get the expressions
 \begin{gather*}
  \alpha_1 =  -
\sum_{i=1}^3 v_i \otimes v_i= - \cos^2 \frac{\gamma}{3}
\sum_{i=1}^3 u_i \otimes u_i - \sin^2 \frac{\gamma}{3}
\sum_{i=1}^3 Ju_i \otimes Ju_i - \sin \frac{2\gamma}{3}
\sum_{i=1}^3 u_i \vee Ju_i,
\\
 \alpha_2 =  -
\sum_{i=1}^3 Jv_i \otimes Jv_i= - \sin^2 \frac{\gamma}{3}
\sum_{i=1}^3 u_i \otimes u_i - \cos^2 \frac{\gamma}{3}
\sum_{i=1}^3 Ju_i \otimes Ju_i  + \sin \frac{2\gamma}{3}
\sum_{i=1}^3 u_i \vee Ju_i,
\end{gather*}
where $\vee$ denotes the symmetric product $a \vee b = 1/2 \, ( a
\otimes b + b \otimes a)$. From equation \eqref{rraa1}, we obtain
\begin{eqnarray*}
r(\nabla^6 \omega) & = & \textstyle - \frac{1}{2} ( \cos \gamma +
\sigma \sin \gamma) \langle \cdot , \cdot \rangle - \frac{1}{2} (
\sin \gamma - \sigma \cos \gamma  ) \omega - ( \sin
\frac{\gamma}{3} + \sigma \cos \frac{\gamma}{3} ) \sum_{i=1}^3 u_i
\vee Ju_i
\\
&& + \textstyle \frac{1}{2} ( \cos  \frac{\gamma}{3} -  \sigma
\sin \frac{\gamma}{3} ) \sum_{i=1}^3 (u_i \otimes u_i - Ju_i
\otimes Ju_i ).
\end{eqnarray*}
The first two terms constitute the $\mathcal W_1$-part of the
tensor $r(\nabla^6 \omega)$, while  the $\mathcal W_3$-part
consists of the last two remaining terms.

We have already deduced at the end of Subsection
\ref{su3structure} that $r( \nabla^6 \omega) = \sum_{j,k=1}^6
c_{jk} e_j \otimes e_k$ implies $ \nabla^6 \omega = \sum_{j,k=1}^6
c_{jk} e_j \otimes e_k \lrcorner \Psi_+$. Then it follows
\begin{gather} \label{nablaexpression}
\hspace{-4cm} \nabla^6 \omega  =  \textstyle - \frac{1}{2}  (\cos
\gamma + \sigma \sin \gamma)  \Psi_+  - \frac{1}{2}  (\sin \gamma
- \sigma \cos \gamma)
\Psi_- \\
\hspace{1.15cm} \textstyle + \frac{1}{2}  (\cos \frac{\gamma}{3} -
\sigma \sin \frac{\gamma}{3})( u_1 \wedge u_2 \wedge u_3 +
\sumcic_{ijk=123} ( u_i \wedge J u_j \wedge J u_k  - 2 u_i \otimes
Ju_j \wedge Ju_k) ) \nonumber\\
\hspace{1.5cm} \textstyle + \frac{1}{2}  (\sin \frac{\gamma}{3} +
\sigma \cos \frac{\gamma}{3}) ( Ju_1 \wedge Ju_2 \wedge Ju_3 +
\sumcic_{ijk=123} ( Ju_i \wedge  u_j \wedge  u_k  - 2 Ju_i \otimes
u_j \wedge u_k) ), \nonumber
\end{gather}
where $\sumcic$ denotes cyclic sum. Thus the exterior derivative
$d \omega$ of the K{\"a}hler form is given by
\begin{eqnarray*}
d \omega & = & \textstyle - \frac{3}{2}  (\cos \gamma + \sigma
\sin \gamma) \Psi_+ - \frac{3}{2}  (\sin \gamma - \sigma \cos
\gamma)
\Psi_- \\
&& \textstyle + \frac{1}{2}  (\cos \frac{\gamma}{3} - \sigma \sin
\frac{\gamma}{3})( 3 u_1 \wedge u_2 \wedge u_3 +
\sumcic_{ijk=123} u_i \wedge J u_j \wedge J u_k   ) \nonumber\\
&& \textstyle + \frac{1}{2}  (\sin \frac{\gamma}{3} + \sigma \cos
\frac{\gamma}{3})( 3 Ju_1 \wedge Ju_2 \wedge Ju_3 +
\sumcic_{ijk=123}  Ju_i \wedge  u_j \wedge  u_k ). \nonumber
\end{eqnarray*}

It was shown in \cite{Cabrera:special} that if $(\nabla^6
\omega)_{\mathcal W_1} = \lambda \Psi_+ + \mu \Psi_-$, then
$\left( d \Psi_+\right)_{\mathcal W_1^a} = 2 \mu \omega \wedge
\omega$ and $\left( d \Psi_-\right)_{\mathcal W_1^a} = 2 \lambda
\omega \wedge \omega$. Combining this with \eqref{dcuatrocinco},
we get from \eqref{nablaexpression} that
\begin{gather} \label{difpsi+}
d \Psi_+ = -  (\sin \gamma - \sigma \cos \gamma) \omega \wedge
\omega + J d \gamma \wedge \Psi_+,\\
d \Psi_- = -  (\cos \gamma + \sigma \sin \gamma) \omega \wedge
\omega  + J d \gamma \wedge \Psi_-. \label{difpsi-}
\end{gather}

In particular, one can consider
$$
\begin{array}{ccl}
 v_1 = - \sigma x_0 e + \sigma x e_0 + x_3 e_1 - x_1 e_3, & \qquad
 &  Jv_1 = \sigma \left( x_6 e_2 + x_5 e_4 - x_4 e_5 - x_2 e_6 \right),  \\
  v_2 = - \sigma x_1 e - x_3 e_0 + \sigma x e_1 + x_0 e_3, & \qquad &
  J v_2 =  \sigma \left( x_4 e_2 - x_2 e_4 + x_6 e_5 - x_5 e_6 \right),   \\
   v_3 = - x_3 e + \sigma  x_1 e_0 - \sigma  x_0 e_1 +  x e_3, &
   \qquad & J v_3 =  - x_5 e_2 + x_6 e_4 + x_2 e_5 - x_4
   e_6.
\end{array}
$$
It is straightforward to check that $\Phi (N_1, v_1, v_2 ,v_3)=1$
and
\begin{gather}\label{mnnov2}
 d v_1 = -2 \sigma  \, v_2 \wedge  v_3, \qquad d v_2 = -2 \sigma  \, v_3 \wedge
v_1, \qquad d v_3 = -2 \sigma \, v_1 \wedge  v_2, \\\nonumber
 d \left( Jv_1 \right) = 2  Jv_2 \wedge  Jv_3, \qquad d \left(Jv_2\right) =2  Jv_3 \wedge
Jv_1, \qquad d \left( Jv_3 \right) = 2  Jv_1 \wedge  Jv_2.
\end{gather}
Now, using \eqref{mnnov1} and \eqref{mnnov2}, we can compute $d
u_i$, $d\left(J u_i \right)$. From these, $d \omega$, $d \Psi_+$
and $d \Psi_-$ can be again obtained
 by an alternative way.

For the Nijenhuis tensor $N$, we calculate
$N= 2\sqrt{2}\, \Psi_-^{\frac{\pi}{4}}$,
where $\Psi_-^{\frac{\pi}{4}}$ is obtained  from \eqref{defsu} for
$\gamma=\frac{\pi}{4}$. Thus, taking $\sigma=+1$ in \eqref{1} and
$\gamma=\frac{\pi}{4}$ in \eqref{defsu},  from \eqref{difpsi+} and
 \eqref{difpsi-}  we obtain
$$
d\Psi_+^{\frac{\pi}{4}} =0, \qquad d\Psi_-^{\frac{\pi}{4}}=
-\sqrt{2}\, \omega\wedge\omega= \textstyle
-\frac18(N,\Psi_-^{\frac{\pi}{4}}) \, \omega\wedge\omega.
$$
 Applying \cite[Theorem 4.1]{II}, we
conclude that the unique $\Lie{U}(3)$-connection
$\widetilde{\nabla}$ with totally skew-symmetric torsion, defined
in \cite{FI}, preserves the $SU(3)$-structure
($\widetilde{\nabla}\Psi_+^{\frac{\pi}{4}}=\widetilde{\nabla}\Psi_-^{\frac{\pi}{4}}=0$)
on $S^3\times S^3$ obtained for $\sigma=+1, \gamma=\frac{\pi}{4}$.
In particular, the Nijenhuis tensor $N$ is
$\widetilde{\nabla}$-parallel and nowhere vanishing. Therefore,
the structure is strict quasi-integrable $\Lie{U}(3)$-structure in
the sense of \cite{Br1}.

More precisely, we have:
\begin{itemize}
\item
 if $\sigma = +1$ and $\gamma =\frac{\pi}{4}, -\frac{3\pi}{4}$, then the $\Lie{SU}(3)$-structure on $ S^3_{\;1}
\times S^3_{\;2}$ is compatible with the standard product metric
and  half-flat of type $\mathcal W^-_1 \oplus \mathcal W_3 $.
 \item
 if $\sigma = -1$ and $\gamma =\frac{\pi}{4}, -\frac{3\pi}{4}$, then the $\Lie{SU}(3)$-structure on $ S^3_{\;1}
\times S^3_{\;2}$ is  compatible with the standard product metric
and half-flat of type $\mathcal W^+_1 \oplus \mathcal W_3 $.
\end{itemize}

Since $S^3_{\;1} \times S^3_{\;2}\subset \mathbb R^8$ is neither
totally umbilic nor minimal, these structures are neither nearly
K\"ahler nor complex. Moreover, for these cases, we have a global
half-flat $\SU(3)$-structure on $S^3_{\;1} \times S^3_{\;2}$ with
totally skew-symmetric $\widetilde{\nabla}$-parallel nowhere
vanishing Nijenhuis tensor. Therefore, each one of such structures
is strict quasi-integrable $\Lie{U}(3)$-structure in the sense of
\cite{Br1} on $S^3\times S^3$ which is neither nearly K\"ahler nor
complex.

\begin{rem}
Consider $S^3\times S^3\cong SU(2)\times SU(2)$ as the group
manifold $SU(2)\times SU(2)$ and observe that the basis defined by
\eqref{mnnov2} is (up to an orientation) the standard
left-invariant basis on the group manifold $SU(2)\times SU(2)\cong
S^3\times S^3$. This shows  that the $U(3)$-structure defined in
Example~\ref{exam1} is left-invariant compatible with the
bi-invariant Riemannian metric on the group $SU(2)\times SU(2)$.
The torsion connection $\widetilde{\nabla}$ coincides with the
flat canonical connection $\bar\nabla$ on the group manifold
$SU(2)\times SU(2)$ defined by making the standard left invariant
basis  $\bar{\nabla}$- parallel.
\end{rem}
\end{example}

\begin{example} The following examples are  already well known, but we pointed
out them just to illustrate results here exposed. We consider the
product manifold of spheres $S^7 \times S^1$. In \cite{C1}, it is
shown that $S^7 \times S^1$ has a locally conformal parallel
$\Spin(7)$-structure such that the  Lee form $\theta^8$ is  a
constant multiple of the Maurer-Cartan one-form on $S^1$. Since
$S^5 \times S^1$ is  a totally geodesic submanifold of $S^7 \times
S^1$ and $\theta^8$ is tangent to $S^5 \times S^1$, by Corollary
\ref{wtwogeneraltan}, the induced $\Un(3)$-structure on $S^5
\times S^1$ is locally conformal K{\"a}hler. On the other hand,
the sphere $S^6$ is totally geodesic in $S^7 \times S^1$, but now
$\theta^8$ is normal to $S^6$. Hence, by Theorem
\ref{wtwogeneral}, the induced $\Un(3)$-structure on $S^6$ is
nearly K{\"a}hler.
\end{example}

\begin{example}
Let  $Hel^2$ be the two-dimensional helicoid
\begin{gather*}
x^0=\sinh u\cos v, \quad x^1=\sinh u\sin v, \quad x^3=v
\end{gather*}
lying in the  Cayley plane  $\mathbb R^4=span\{e,e_0,e_1,e_3\}$.
 Taking the frame on the normal bundle
$$
 N_1 \cosh u= - \sin v e_0 + \cos v e_1 - \sinh u e_3, \qquad N_2
= e,
$$
 the $SU(3)$-structure on $M^6=Hel^2\times Q^4$ induced by the
standard $\Spin(7)$-structure \eqref{1} on $\mathbb R^8$, $Q^4 =
span\{e_2,e_4,e_5,e_6\}$, is given by the equations
\begin{eqnarray*}
\omega \cosh u& = &  \cos v(e_2\wedge e_4+e_5\wedge e_6)
- \sin v(e_2\wedge e_6+e_4\wedge e_5)\\
&&- \sinh v(e_4\wedge e_6-e_2\wedge e_5)-\cosh^3 u \, du\wedge dv,\\
 \Psi_+=  N_1 \lrcorner \Phi & = & -(-\sinh u\cos v du + \cosh
u\sin v
dv)\wedge(e_2\wedge e_4+e_5\wedge e_6)\\
 && - (\sinh u\sin v du +
\cosh u\cos v dv)\wedge(e_2\wedge
e_6+e_4\wedge e_5)\\
& &- du\wedge(e_2\wedge e_5-e_4\wedge e_6), \\
 \Psi_- =  \sigma N_2
\lrcorner \Phi & = & (\cosh u\sin v du + \sinh u\cos
v dv)\wedge(e_2\wedge e_4+e_5\wedge e_6) \\
 && + (\cosh u\cos v du -
\sinh u\sin v dv)\wedge(e_2\wedge e_6+e_4\wedge e_5) \\
&& - dv\wedge(e_2\wedge e_5-e_4\wedge e_6).
\end{eqnarray*}
Clearly this structure is holomorphic, $d\Psi_\pm=0$ with zero Lee
form, $\theta^6=0$. Therefore, the Bismut connection preserves
this $\SU(3)$-structure due to \cite[Theorem 4.1]{II}, i.e. it has
holonomy contained in $\SU(3)$.

We note that if the helicoid does not lie in a Cayley plane the
induced $\SU(3)$-structure could be not closed.
\end{example}

\begin{table}[tbp]
  \centering {\footnotesize
  \begin{tabular}{ll}
    \toprule
    $ {\mathcal W}^+_1 + {\mathcal W}_1^- + {\mathcal W}^+_2 + {\mathcal W}_2^- +
{\mathcal W}_3 $ &  $J d \gamma=  \frac12 \ast_6 \left( \ast_6 f^*(
L_{N_1} \Phi) \wedge f^* (N_1 \lrcorner \Phi) \right)$
         \\
    \midrule
    ${\mathcal W}^+_1 + {\mathcal W}_1^- + {\mathcal W}^+_2 + {\mathcal W}_2^-
    + {\mathcal W}_5$  &  $(1-J) \sigma \alpha_1 = J_{(1)}  \left( 1 -J \right) \alpha_2$\\
     \midrule
     ${\mathcal W}^+_1 + {\mathcal W}_1^- + {\mathcal W}^-_2 +  {\mathcal W}_3 + {\mathcal W}_5$  &
     $ \cos \gamma (1+J) \sigma \alpha_1 - \sin \gamma \, (1+J) \alpha_2 =
       2\left(\sigma h_1 \cos \gamma -  h_2 \sin \gamma \right) \langle \cdot , \cdot \rangle$ \\
      \midrule
      ${\mathcal W}^+_1 + {\mathcal W}_1^-
      + {\mathcal W}^+_2 +  {\mathcal W}_3 + {\mathcal W}_5$ &
 $ \sin \gamma (1+J) \sigma \alpha_1 + \cos \gamma \, (1+J) \alpha_2
 = 2\left(\sigma h_1 \sin \gamma +  h_2 \cos \gamma \right) \langle \cdot , \cdot \rangle$
      \\
      \midrule
      ${\mathcal W}^-_1 + {\mathcal W}_2^+ +
       {\mathcal W}^-_2 +  {\mathcal W}_3 + {\mathcal W}_5$ &   $\sigma
       h_1 \cos \gamma= h_2 \sin \gamma$ \\
      \midrule
      ${\mathcal W}^+_1 + {\mathcal W}_2^+ + {\mathcal W}^-_2 +  {\mathcal W}_3 + {\mathcal W}_5$
      & $\sigma       h_1 \sin \gamma= - h_2 \cos \gamma$ \\
      \midrule
      ${\mathcal W}_1^+ + {\mathcal W}^-_1 +  {\mathcal W}_3 + {\mathcal W}_5$ &
      $(1+J) \alpha_1  =     2h_1      \langle \cdot , \cdot \rangle$ and
      $(1+J)  \alpha_2 =       2h_2      \langle \cdot , \cdot \rangle$ \\
      \midrule
       ${\mathcal W}_2^+ + {\mathcal W}^-_2 + {\mathcal W}_3 + {\mathcal W}_5$ & $h_1=0$ and $h_2=0$, i.e.
       $M^6$ is  minimal  \\
       \midrule
      ${\mathcal W}^+_1+ {\mathcal W}^-_1 + {\mathcal W}_5$ & $\alpha_1
      = h_1      \langle \cdot , \cdot \rangle$ and $\alpha_2 = h_2      \langle \cdot , \cdot \rangle$, i.e. $M^6$ is totally umbilic  \\
       \midrule
       ${\mathcal W}_3 + {\mathcal W}_5$ & $J \alpha_1=-\alpha_1$ and  $J \alpha_2=-\alpha_2$, in particular, $M^6$ is  minimal \\
       \midrule
       ${\mathcal W}_5$ &  $M^6$ is totally geodesic \\
    \bottomrule
  \end{tabular}  } \vspace{2mm}
  \caption{$M^8$ of type parallel ($\overline{W}_0$)}
  \label{tab:wzerogeneral}
\end{table}

\begin{table}[tbp]
  \centering {\footnotesize
  \begin{tabular}{ll}
    \toprule
 ${\mathcal W}^+_1 + {\mathcal W}_1^- + {\mathcal W}^+_2 + {\mathcal
 W}_2^-+ {\mathcal W}_3
 + {\mathcal W}_5$  & $\overline{r} ( \nabla^8\Phi )(N_1,N_1, f_{\ast} \cdot ) = -
\sigma \overline{r} (\nabla^8 \Phi )(N_2, f_{\ast}  J \cdot, N_1
)$ \\
& \hspace{3.1cm} $  +\sigma \overline{r} (\nabla^8 \Phi )(f_{\ast}  J \cdot, N_2, N_1 )$ \\
     \midrule
      ${\mathcal W}^-_1 + {\mathcal W}_2^+ +
       {\mathcal W}^-_2 +  {\mathcal W}_3 + {\mathcal W}_4 + {\mathcal W}_5$ &
 $ \sin \gamma \left( \sigma h_1 - \sigma
\overline{r} (\nabla^8 \Phi )(N_2,N_2,N_1) \right)  =$ \\
& \hspace{3.2cm} $ =\cos \gamma \left(
 h_2 -  \overline{r} (\nabla^8 \Phi )(N_1,N_2,N_1) \right)$
       \\
      \midrule
      ${\mathcal W}^+_1 + {\mathcal W}_2^+ + {\mathcal W}^-_2 +  {\mathcal W}_3 + {\mathcal W}_4  + {\mathcal W}_5$
      & $  \cos \gamma \left( \sigma h_1 - \sigma
\overline{r} (\nabla^8 \Phi )(N_2,N_2,N_1) \right)  = $ \\
& \hspace{3.2cm}  $= - \sin \gamma \left(
 h_2 -  \overline{r} (\nabla^8  \Phi )(N_1,N_2,N_1) \right)$ \\
      \midrule
      ${\mathcal W}_2^+ + {\mathcal W}^-_2 +  {\mathcal W}_3 +  {\mathcal W}_4  + {\mathcal W}_5$ &
       $h_1 = \overline{r} (\nabla^8  \Phi )(N_2,N_2,N_1)$ and $h_2 =  \overline{r} (\nabla^8  \Phi )(N_1,N_2,N_1)$  \\
    \bottomrule
  \end{tabular}  } \vspace{2mm}
  \caption{$M^8$ of type balanced ($\overline{W}_1$)}
  \label{tab:wonegeneral}
\end{table}

\begin{table}[tbp]
  \centering {\footnotesize
  \begin{tabular}{ll}
    \toprule
    $ {\mathcal W}^+_1 + {\mathcal W}_1^- + {\mathcal W}^+_2 + {\mathcal W}_2^- + {\mathcal W}_3 + {\mathcal W}_4$ &
      $2J d \gamma=   \ast_6 \left( \ast_6 f^*(
L_{N_1} \Phi) \wedge f^* (N_1 \lrcorner \Phi) \right) + \theta^6$
         \\
    \midrule
    ${\mathcal W}^+_1 + {\mathcal W}_1^- + {\mathcal W}^+_2 + {\mathcal W}_2^-
    + {\mathcal W}_3 + {\mathcal W}_5$  & $\theta^8$ is normal to $M^6$ \\
    \midrule
    ${\mathcal W}^+_1 + {\mathcal W}_1^- + {\mathcal W}^+_2 + {\mathcal W}_2^-
    + {\mathcal W}_4 + {\mathcal W}_5$  &
    $\sigma(1-J) \alpha_1 = J_{(1)}  \left( 1 -J \right) \alpha_2$\\
     \midrule
     ${\mathcal W}^+_1 + {\mathcal W}_1^- + {\mathcal W}^-_2 +  {\mathcal W}_3 + {\mathcal W}_4 + {\mathcal W}_5$  &
     $ \cos \gamma (1+J) \sigma \alpha_1 - \sin \gamma \, (1+J) \alpha_2  =
         2\left(\sigma h_1 \cos \gamma -  h_2 \sin \gamma \right) \langle \cdot , \cdot \rangle$ \\
      \midrule
      ${\mathcal W}^+_1 + {\mathcal W}_1^-
      + {\mathcal W}^+_2 +  {\mathcal W}_3  + {\mathcal W}_4 + {\mathcal W}_5$ &
 $ \sin \gamma (1+J) \sigma \alpha_1 + \cos \gamma \, (1+J) \alpha_2
 = 2\left(\sigma h_1 \sin \gamma +  h_2 \cos \gamma \right) \langle \cdot , \cdot \rangle$
      \\
      \midrule
      ${\mathcal W}^-_1 + {\mathcal W}_2^+ +
       {\mathcal W}^-_2 +  {\mathcal W}_3  + {\mathcal W}_4 + {\mathcal W}_5$ &   $\sigma \left(
       \theta^8(N_1) - 4  h_1 \right) \cos \gamma= \left(\theta^8(N_2) - 4 h_2 \right) \sin \gamma$ \\
      \midrule
      ${\mathcal W}^+_1 + {\mathcal W}_2^+ + {\mathcal W}^-_2 +  {\mathcal W}_3  + {\mathcal W}_4 + {\mathcal W}_5$
      & $\sigma \left(  \theta^8(N_1) - 4 h_1 \right) \sin \gamma= - \left(  \theta^8(N_2) - 4h_2 \right) \cos \gamma$ \\
      \midrule
      ${\mathcal W}_1^+ + {\mathcal W}^-_1 +  {\mathcal W}_3  + {\mathcal W}_4 + {\mathcal W}_5$ &
       $(1+J) \alpha_1  =   2h_1      \langle \cdot , \cdot \rangle$
       and $(1+J)  \alpha_2 =  2h_2      \langle \cdot , \cdot \rangle$ \\
      \midrule
       ${\mathcal W}_2^+ + {\mathcal W}^-_2 + {\mathcal W}_3  + {\mathcal W}_4 +
{\mathcal W}_5$ & $4h_1=\theta^8(N_1)$ and $4h_2= \theta^8(N_2)$
              \\
      \midrule
          ${\mathcal W}^+_1+ {\mathcal W}^-_1 + {\mathcal W}_4 + {\mathcal W}_5$
          &$M^6$ is  totally umbilic  \\
\midrule
 $ {\mathcal W}_4 + {\mathcal W}_5$ &
$4\alpha_1
      =  \theta^8(N_1)      \langle \cdot , \cdot \rangle$ and  $4\alpha_2 = \theta^8(N_2)      \langle \cdot , \cdot \rangle$
\\
\midrule $ {\mathcal W}_5$ &
$4\alpha_1
      = \theta^8(N_1)      \langle \cdot , \cdot \rangle$,  $4\alpha_2 =  \theta^8(N_2)      \langle \cdot , \cdot
     \rangle$ and  $\theta^8$ is normal to $M^6$ \\
    \bottomrule
  \end{tabular}  } \vspace{2mm}
  \caption{$M^8$ of type locally conformal parallel ($\overline{W}_2$)}
  \label{tab:wtwogeneral}
\end{table}

\begin{table}[tbp]
  \centering {\footnotesize
  \begin{tabular}{ll}
    \toprule
    $ {\mathcal W}^+_1 + {\mathcal W}_1^- + {\mathcal W}^+_2 + {\mathcal W}_2^- + {\mathcal W}_3 + {\mathcal W}_4$ &
      $2J d \gamma=   \ast_6 \left( \ast_6 f^*(
L_{N_1} \Phi) \wedge f^* (N_1 \lrcorner \Phi) \right) + \theta^6$
         \\
    \midrule
    ${\mathcal W}^+_1 + {\mathcal W}_1^- + {\mathcal W}^+_2 + {\mathcal W}_2^-
    + {\mathcal W}_3 + {\mathcal W}_5$  & $f^* \theta^8 =0$  \\
    \midrule
    ${\mathcal W}^+_1 + {\mathcal W}_1^- + {\mathcal W}^+_2 + {\mathcal W}_2^-
    + {\mathcal W}_4 + {\mathcal W}_5$  &
    $\sigma(1-J) \alpha_1 = J_{(1)}  \left( 1 -J \right) \alpha_2$\\
     \midrule
     ${\mathcal W}^+_1 + {\mathcal W}_1^- + {\mathcal W}^-_2 +  {\mathcal W}_3 + {\mathcal W}_4 + {\mathcal W}_5$  &
     $ \cos \gamma (1+J) \sigma \alpha_1 - \sin \gamma \, (1+J) \alpha_2  =     2\left(\sigma h_1 \cos \gamma -  h_2 \sin \gamma \right) \langle \cdot , \cdot \rangle$ \\
      \midrule
      ${\mathcal W}^+_1 + {\mathcal W}_1^-
      + {\mathcal W}^+_2 +  {\mathcal W}_3  + {\mathcal W}_4 + {\mathcal W}_5$ &
 $ \sin \gamma (1+J) \sigma \alpha_1 + \cos \gamma \, (1+J) \alpha_2
 = 2\left(\sigma h_1 \sin \gamma +  h_2 \cos \gamma \right) \langle \cdot , \cdot \rangle$
      \\
      \midrule
      ${\mathcal W}^-_1 + {\mathcal W}_2^+ +
       {\mathcal W}^-_2 +  {\mathcal W}_3  + {\mathcal W}_4 + {\mathcal W}_5$ &   $\sigma  h_1 \cos \gamma
       = h_2  \sin \gamma$ \\
      \midrule
      ${\mathcal W}^+_1 + {\mathcal W}_2^+ + {\mathcal W}^-_2 +  {\mathcal W}_3  + {\mathcal W}_4 + {\mathcal W}_5$
      & $\sigma h_1  \sin \gamma= - h_2  \cos \gamma$ \\
      \midrule
      ${\mathcal W}_1^+ + {\mathcal W}^-_1 +  {\mathcal W}_3  + {\mathcal W}_4 + {\mathcal W}_5$ &
      $(1+J) \alpha_1       = 2h_1      \langle \cdot , \cdot \rangle$ and $(1+J)  \alpha_2  =
        2h_2      \langle \cdot , \cdot \rangle$ \\
      \midrule
       ${\mathcal W}_2^+ + {\mathcal W}^-_2 + {\mathcal W}_3  + {\mathcal W}_4 +
{\mathcal W}_5$ & $M^6$ is minimal
              \\
      \midrule
          ${\mathcal W}^+_1+ {\mathcal W}^-_1 + {\mathcal W}_4 + {\mathcal W}_5$
          &$M^6$ is  totally umbilic  \\
\midrule
 $ {\mathcal W}_4 + {\mathcal W}_5$ &
$M^6$ is totally geodesic
\\
\midrule $ {\mathcal W}_5$ &
$M^6$ is totally geodesic and   $f^*\theta^8=0$ \\
    \bottomrule
  \end{tabular}  } \vspace{2mm}
  \caption{$M^8$ of type locally conformal parallel ($\overline{W}_2$) with Lee form $\theta^8$ tangent to $M^6$}
  \label{tab:wtwogeneraltan}
\end{table}

\bibliographystyle{hamsplain}

\providecommand{\bysame}{\leavevmode\hbox
to3em{\hrulefill}\thinspace}






\end{document}